\documentclass[11pt]{article}
\usepackage{amsmath}
\usepackage[table]{xcolor}
\usepackage{graphicx}
\usepackage{amsfonts}
\usepackage{amsthm}
\usepackage{tikz}
\usepackage{epsfig,epsf,psfrag}

\usepackage{subcaption}
\usepackage{caption}

\newcommand{\double}{\baselineskip 1.24 \baselineskip}
\newcommand\valpha{\vec \alpha}
\newcommand\vsigma{\vec \sigma}
\newcommand\vbeta{\vec \beta}
\newcommand\veca{\vec a}

\newtheorem{theorem}{Theorem}[section]
\newtheorem{lemma}[theorem]{Lemma}

\theoremstyle{definition}
\newtheorem{definition}{Definition}[section]

\theoremstyle{remark}
\newtheorem{remark}{Remark}

\def\keywords{\vspace{.5em}
{\textit{Keywords}:\,\relax%
}}

\title{A fast solver for multi-particle scattering in a layered medium}
\author{\begin{tabular}{c}
Jun Lai\thanks{Courant Institute of Mathematical Sciences, New York University, NY 10012 
(Email: lai@cims.nyu.edu, motokia.kobayashi@gmail.com, greengar@cims.nyu.edu) } \and Motoki Kobayashi\footnotemark[1] \and Leslie Greengard\footnotemark[1] \thanks{Simons Center for Data Analysis, Simons Foundation, New York, NY 10010}
\end{tabular}}
\date{}

\begin{document}
\maketitle \double
\abstract{ In this paper, we consider acoustic or electromagnetic scattering in two dimensions from an infinite 
three-layer medium with thousands of wavelength-size dielectric particles embedded in the middle layer. 
Such geometries are typical of microstructured composite materials, and the evaluation of the scattered
field requires a suitable fast solver for either a single configuration or for a sequence of 
configurations as part of a design or optimization process. 
We have developed an algorithm for problems of this type by combining the Sommerfeld integral representation, 
high order integral equation discretization, the fast multipole method and classical multiple scattering theory. 
The efficiency of the solver is illustrated with several numerical experiments.}

\section{Introduction}

The problem of designing composite materials that exhibit a specific acoustic or
electromagnetic response is an area of
active research \cite{BJL1,Parnell01052010, Wuzhang2009}. 
Examples include the design of random media with a well-defined macroscopic refraction (coherent scattering) 
\cite{Parnell01052010} and the fabrication of metamaterials \cite{Wuzhang2009} for cloaking, near field imaging,
etc.  In many experiments, the materials are designed by incorporating large numbers of identical
inclusions (particles) in a layered material. 
When the size of each particle is comparable to the wavelength of the incoming field and the distribution of 
particles is reasonably dense, then
the interaction of the particles involves non-negligible multiple scattering effects and 
methods based on homogenization \cite{Parnell01052010} are not applicable.
Instead, the full Helmholtz or Maxwell equations should be solved at each iteration of the design
process.  Numerical simulation, in the absence of suitable fast algorithms, are impractical when
thousands of particles are involved. 

In this paper, we develop an algorithm that accelerates the computation of electromagnetic scattering when a 
large number of particles are embedded in the middle of a three-layer dielectric medium. 
Numerical experiments show that our solver takes 1--2 minutes to evaluate the scattered field for up to 
$5,000$ particles on a $2.3$GHz laptop. Our method combines the Sommerfeld integral representation, 
a well-posed integral formulation, high-order discretization, multiple scattering theory and the fast
multipole method. We focus on the two dimensional setting by assuming the material is invariant in the $z$ 
direction. A related three-dimensional solver was considered in \cite{GG2013}, but
the particles were assumed to be distributed in free space. A principal contribution of this paper
is the development of a mathematical framework that permits them to be embedded in a layer
material (which is closer to being manufacturable).
While we restrict our attention here to the three-layer case, the extension to an arbitrary
layered medium is straightforward. 

\begin{figure}[tbp]
\begin{center}
\includegraphics[width=80mm]{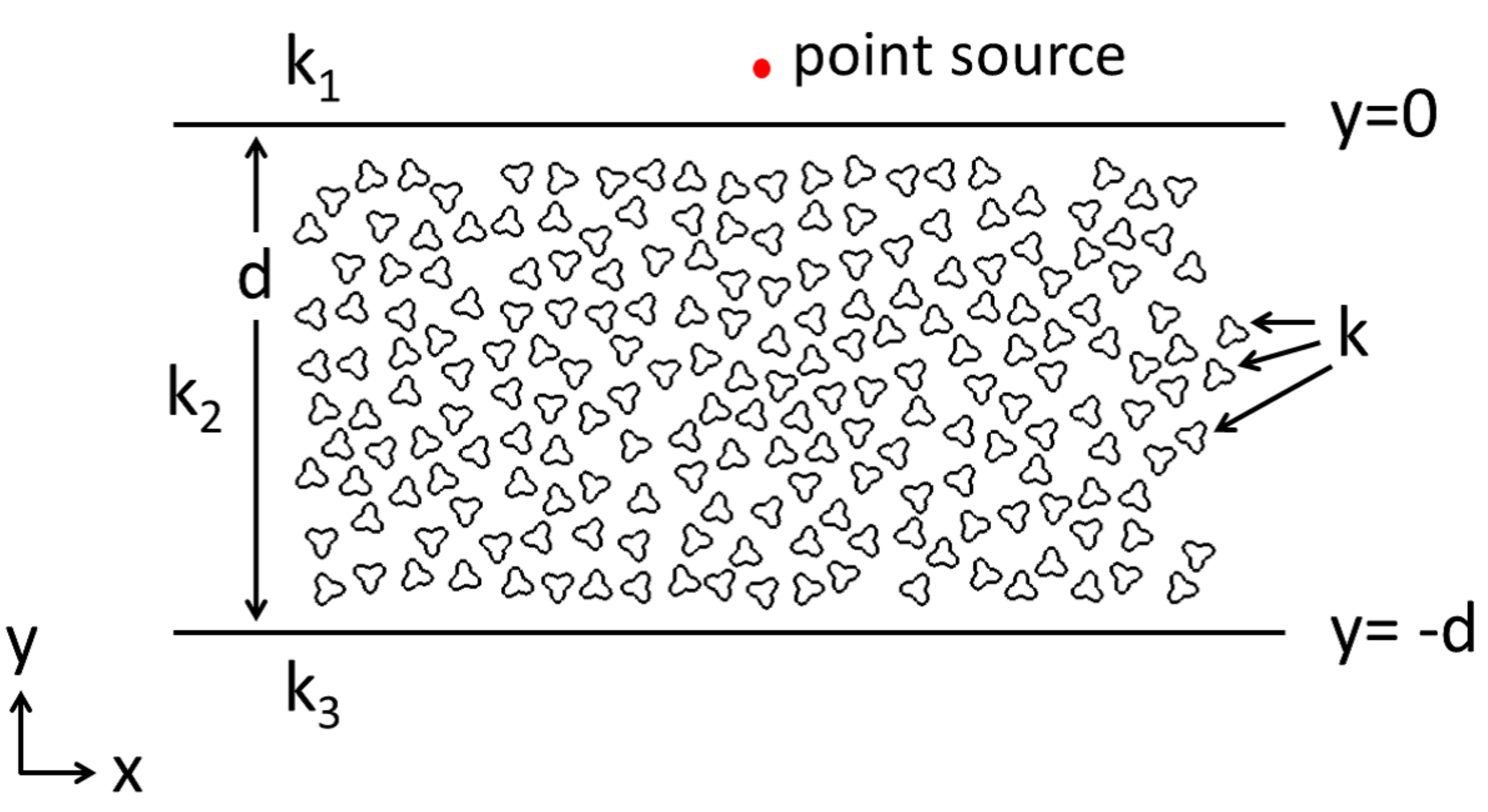}
\caption{Geometry of the three-layered medium, with a large number of dielectric particles embedded in the 
middle layer.}
\label{fig:NumericalModel}
\end{center}
\end{figure}

More precisely, we consider time-harmonic scattering (with time dependence $e^{i\omega t}$) from
a three-layered medium as depicted in Fig. \ref{fig:NumericalModel}. The incident field is 
assumed to be driven by a point source located in the first (top) layer. 
The thickness of the middle layer is denoted by $d$. We assume the magnetic permeability $\mu$ is identical 
in each layer, while the electric permittivity $\epsilon$ is piecewise constant. 
There are two 
fundamental polarizations in the two dimensional setting to consider: 
the transverse magnetic (TM) polarization and 
the transverse electric (TE) polarization. In both cases, the Maxwell equations reduce to a scalar
Helmholtz equation. For simplicity, we consider the TM polarization here, in which case the scattered field 
$u^s$ must satisfy the equation

\begin{equation}
\label{eqintr}
\Delta u^s+k^2 u^s = 0,
\end{equation} 
where $k = \omega\sqrt{\mu\epsilon}$ is the wavenumber. We 
denote by $k_1,\ k_2$ and $k_3$ the wavenumber for the three layers, and by $k_p$ 
the wavenumber for the particles. 
The scattered field also has to satisfy the Sommerfeld radiation condition at infinity \cite{Cot2}:
\begin{equation}\label{SommerRadCond}
\lim_{r\rightarrow\infty} \sqrt{r} \bigg(\frac{\partial u^s}{\partial n}-iku^s\bigg) = 0,
\end{equation}
where $r=\sqrt{x^2+y^2}$.

In order to develop an especially fast solution method, we make two further 
assumptions.
First, as in the fast multi-particle scattering (FMPS) method of 
\cite{GG2013},  we assume that the particles are well separated 
from each other - that is, the 
separation between particles is at least 10\% of the particle size. Second,
we assume that only a finite number of distinct particle shapes 
are included in the simulation. The first condition
ensures that a multiple scattering formalism will be accurate and the second 
condition ensures that 
precomputation of single particle scattering matrices permits a
dramatic reduction in the number of degrees of freedom necessary for the solver. 
The particles are not assumed to be symmetric and may be placed with
arbitrary orientation. Both hypotheses are common in materials design 
(although there are exceptions).

An outline of the paper follows. In Section 2,  we introduce the Sommerfeld integral and its application
to layered materials (in the absence of inclusions). In Section 3, we review classical multiple scattering theory 
for circular particles. Section 4 extends the scattering formalism to non-circular 
particles and Section 5 
develops analytical tools needed to go back and forth between the Sommerfeld integral formalism and multiple
scattering theory.  Section 5 also combines the techniques in the preceding sections 
and extends the FMPS
method to layered media. 
Numerical examples are provided in Section 6 to 
illustrate the efficiency of the method, followed by some concluding remarks in 
Section 7.

\section{The Sommerfeld integral for layered media}
\label{sec2}
Wave propagation in a layered medium is a well-studied problem in 
acoustic and electromagnetic scattering theory. 
Nearly a century ago, Sommerfeld developed a spectral representation involving a
Fourier integral in the ``transverse" variable (the $x$-coordinate in Fig.
\ref{fig:NumericalModel}) \cite{CWC1995}. 
Assuming a point source is located at $\mathbf{x_0}=(x_0,y_0)$ in the top layer,
 with wavenumber $k_1$, the corresponding field is given by the 
(two-dimensional) free space Green's function: 
$G_{k_1}(\mathbf{x},\mathbf{x_0}) = 
\frac{i}{4}H_0^{(1)}(k_1|\mathbf{x}-\mathbf{x_0}|)$, 
where $H_0^{(1)}(x)$ is the first kind Hankel function of order zero. 
Combing the Fourier transform and contour integration \cite{ONeil2014}, 
the Green's function can also be written in the form:
\begin{align}\label{sommer}
G_{k_1}(\mathbf{x},\mathbf{x_0}) = \frac{1}{4\pi}\int_{-\infty}^{\infty}\frac{e^{-\sqrt{\lambda^2-{k_1}^2}|y-y_0|}}{\sqrt{\lambda^2-{k_1}^2}} e^{i\lambda(x-x_0)}d\lambda.
\end{align}
It is important to note that the Sommerfeld integral (\ref{sommer}) 
is conditionally convergent and as stated, requires that
$y\neq y_0$.

In the Sommerfeld approach (\cite{CWC1995}),
we assume the upward scattered field $u_1^s$ in the top layer can be 
expressed as
\begin{equation}\label{layer1}
u_1^s = \frac{1}{4\pi}\int_{-\infty}^{\infty}\frac{e^{-\sqrt{\lambda^2-{k_1}^2}y}}{\sqrt{\lambda^2-{k_1}^2}} e^{i\lambda(x-x_0)}\sigma_1(\lambda)d\lambda,
\end{equation} 
where $\sigma_1(\lambda)$ is an unknown density on the upper interface $y = 0$.
It is straghtforward to verify that $u_1^s$ satisfies the Helmholtz
equation with Helmholtz parameter $k_1$.

In the second layer, the scattered field $u_2^s$ can be written in terms
of contributions from
both the upper ($y=0$) and lower ($y=-d$) interfaces: $u_2^t$ and $u_2^b$. 
These are given by
\begin{align}
u_2^t &= \frac{1}{4\pi}\int_{-\infty}^{\infty}\frac{e^{\sqrt{\lambda^2-{k_2}^2}y}}{\sqrt{\lambda^2-{k_2}^2}} e^{i\lambda(x-x_0)}\sigma_2^+(\lambda)d\lambda, \label{layer21}\\
u_2^b &= \frac{1}{4\pi}\int_{-\infty}^{\infty}\frac{e^{-\sqrt{\lambda^2-{k_2}^2}(y+d)}}{\sqrt{\lambda^2-{k_2}^2}} e^{i\lambda(x-x_0)}\sigma_2^-(\lambda)d\lambda, \label{layer22}
\end{align}
where $\sigma_2^+(\lambda)$ and 
$\sigma_2^-(\lambda)$ are used to denote {\em spectral density} functions
on the upper and lower interfaces.

Similarly, we can represent the scattered field $u_3^s$ in the third layer with 
an unknown density $\sigma_3(\lambda)$ on the lower interface as
\begin{equation}\label{layer3}
u_3^s = \frac{1}{4\pi}\int_{-\infty}^{\infty}\frac{e^{\sqrt{\lambda^2-{k_3}^2}(y+d)}}{\sqrt{\lambda^2-{k_3}^2}} e^{i\lambda(x-x_0)}\sigma_3(\lambda)d\lambda.
\end{equation}

\begin{remark}
The signs of the terms
$e^{\pm \sqrt{\lambda^2-{k_i}^2}y}$ and $e^{\pm \sqrt{\lambda^2-{k_i}^2}(y+d)}$
in eqs. \eqref{layer1}-\eqref{layer3} 
ensure that evanescent modes 
(when $|\lambda| > |k_i|$) decay away from each layer. 
(Physically, this is related to causality and is 
required in the derivation of formula \cite{ONeil2014} by countour integration).
\end{remark}

It is worth noting that the 
four unknown densities $\sigma_1$, $\sigma_2^+$, $\sigma_2^-$ and $\sigma_3$ can be 
interpreted in two ways. First, they can simply be considered the spectral densities
in the Fourier domain of a consistent representation for the Helmholtz equation.
For those more familiar with potential theory, they can be viewed as the 
Fourier transforms of 
charge densities of four single layer potentials lying on the 
corresponding interfaces \cite{BG2011}. 

In the absence of any inclusions, the Sommerfeld representation for the field in
each subdomain is derived from a ``mode by mode" analysis. That is,
the unknown functions $\sigma_1$, 
$\sigma^+_2$, $\sigma^-_2$, and $\sigma_3$ are found
by enforcing the continuity conditions at the interface for each value of the 
argument $\lambda$. For the case of electromagnetic scattering in TM
polarization, when the permeability $\mu$ is constant in each layer, this 
requires that
\begin{align}
[u] = 0,\label{cont1}\\
\left[\frac{\partial u}{\partial n} \right] = 0,\label{cont2}
\end{align}
where $[\cdot]$ denotes the jump of a function along the interface, 
$\partial /{\partial n}$ is the normal derivative and $u$ is the total field 
in each layer \cite{Cot2}.

It is straightforward to check that
the linear system to be solved for each $\lambda$ takes the form:

\begin{equation}
\left(
\begin{array}{cccc}
\frac{1}{\sqrt{\lambda^2-k_1^2}} & -\frac{1}{\sqrt{\lambda^2-k_2^2}} 
& -\frac{e^{-\sqrt{\lambda^2-k_2^2}d}}{\sqrt{\lambda^2-k_2^2}} & 0 \\
\\
0 & \frac{e^{-\sqrt{\lambda^2-k_2^2}d}}{\sqrt{\lambda^2-k_2^2}} &
\frac{1}{\sqrt{\lambda^2-k_2^2}}  & -\frac{1}{\sqrt{\lambda^2-k_3^2}} \\
\\
1 & 1 & -e^{-\sqrt{\lambda^2-k_2^2}d} & 0 \\
\\
0 & e^{-\sqrt{\lambda^2-k_2^2}d} & -1 & -1 
\end{array} \right)
\left( \begin{array}{c}
\sigma_1(\lambda) \\
\\
\sigma^+_2(\lambda) \\
\\
\sigma^-_2(\lambda) \\
\\
\sigma_3(\lambda) 
\end{array} \right)
=
\left( \begin{array}{cccc}
-\frac{e^{-\sqrt{\lambda^2-k_1^2}y_0}}{\sqrt{\lambda^2-k_1^2}} \\
\\
0 \\
\\
e^{-\sqrt{\lambda^2-k_1^2}y_0} \\
\\
0 
\end{array} \right)
\label{sommerfeldsys}
\end{equation}

\vspace{.2in}

\begin{definition}
We will denote the $4\times 4$ matrix above by $A_\lambda$.
\end{definition}

For the problem we consider here, the Sommerfeld integrals 
must be coupled to a representation of the field induced by the many 
particles present in the central layer. 
Before discussing the coupled system, however, we first summarize some well-known
facts about scattering from a finite collection of inclusions in a homogeneous 
infinite medium. 

\section{Wave scattering for disks}
\label{sec3}
Suppose now that we have an inclusion of dielectric material with 
$k = \omega\sqrt{\epsilon_p \mu}$ embedded in ${\bf R}^2$, assumed to consist of 
a dielectric with 
$k_2 = \omega\sqrt{\epsilon_2 \mu}$.
For transverse magnetic(TM) polarization, the total electrical field $u$ in the 
exterior of the inclusion satisfies the Helmholtz equation:
\begin{equation}\label{ex}
\Delta u + k_2^2u = 0.
\end{equation}   
Further, the total field $u$ can be written as the sum of the incident field 
$u^{inc}$ and the scattered field $u^s$, where $u^s$ satisfies 
\eqref{ex} and the Sommerfeld radiation condition,
\begin{equation}
\lim_{r\rightarrow \infty} \sqrt{r} \bigg(\frac{\partial u^s}{\partial n}-i k_2 u^s\bigg) = 0,
\end{equation}
where $r=\sqrt{x^2+y^2}$.
Within the inclusion, the field $u$ satisfies the Helmholtz equation with 
wavenumber $k_p$,
\begin{equation}\label{in}
\Delta u+k_p^2u = 0.
\end{equation}
On the boundary of the inclusion, we must enforce the continuity conditions 
given by Eq (\ref{cont1}) and (\ref{cont2}).

\subsection{A single disk}

When the inclusion is a disk of radius $R$ centered at the origin, 
it is straightforward 
to represent the solution using separation of variables, with 
\begin{equation}\label{out}
u^s = \sum_{n=-\infty}^{\infty} \beta_n H_n(k_2 r)e^{i n\theta}
\end{equation}
in the exterior and 
\begin{equation}\label{int}
u = \sum_{n=-\infty}^{\infty} \gamma_n J_n(k_p r)e^{i n\theta}
\end{equation}
in the interior. Here, $(r,\theta)$ are the polar coordinates of a point in the
plane, $H_n(r)$ is the Hankel function of the first kind of order $n$ and
$J_n(r)$ is the Bessel function of order $n$ \cite{Cot,Hand2010}.

We now 
expand the incident wave $u^{inc}$ and its normal derivative in the form:
\begin{equation}\label{incident}
u^{inc} = \sum_{n=-\infty}^{\infty} \alpha_n J_n(k_2 r) e^{i n\theta}, \quad
\frac{\partial u^{inc}}{\partial r}  = \sum_{n=-\infty}^{\infty} \alpha_n k_2 J_n'(k_2 r) e^{i n\theta}.
\end{equation}
Enforcing the continuity conditions \eqref{cont1}, \eqref{cont2} on the 
boundary of the disk for each Fourier mode, 
we easily obtain the following linear equation for mode $n$:
\begin{equation}\label{linear}
\left[ \begin{array}{cc}
- H_n(k_2R) &  J_n(k_p R) \\ 
-k_2{H_n}'(k_2R) & k_pJ_n'(k_p R)
\end{array} \right] 
\left[ \begin{array}{c}
\beta_n\\
\gamma_n
\end{array} \right]=
\left[\begin{array}{c}
\alpha_nJ_n(k_2 R) \\
\alpha_nk_2J_n'(k_2 R)
\end{array} \right],
\end{equation}
where $n\in\mathbb{N}$.

Solving Eq. \eqref{linear} determines the coefficients $\beta_n, \gamma_n$: 
\begin{align}
\beta_n &= \left[\frac{k_p J_n(k_2 R)J_n'(k_p R)-k_2J_n'(k_2 R)J_n(k_p R)
}{k_2{H_n}'(k_2R)J_n(k R)-
k_pJ_n'(k_p R)H_n(k_2R)} \right] \alpha_n , \\
\gamma_n & = \left[ \frac{k_2J_n(k_2 R){H_n}'(k_2R)-
k_2J_n'(k_2 R)H_n(k_2R)}{k_2{H_n}'(k_2R)J_n(k_p R)-
k_pJ_n'(k_p R)H_n(k_2R)} \right] \alpha_n.
\end{align}
It is straightforward to verify that the denominator in the preceding expressions
cannot vanish if $k$ and $k_p$ have positive real part and non-negative imaginary part \cite{Cot2,KR78}.

\begin{definition}
The mapping between the incoming coefficients $\{\alpha_n\}$ and outgoing coefficients $\{\beta_n\}$ is referred as the scattering matrix for the disk and denoted by $S$.
\end{definition}

\begin{remark}
While we restrict our attention here to dielectric particles, the method can 
easily be extended to perfectly conducting disks. 
Since the interior field $u$ in \eqref{int} is zero for perfect conductors, 
from \eqref{linear}, we have
\begin{equation}
\beta_n = -\frac{J_n(k_2 R)}{H_n(k_2R)} \, \alpha_n.
\end{equation}
\end{remark}

\begin{remark}
In the remainder of this paper, we will refer expansions based on 
Hankel functions, such as \eqref{out}, as multipole expansions or $H$-expansions, 
and expansions based on Bessel functions, such as \eqref{int}, as local expansions
or $J$-expansions.
\end{remark}	

\begin{remark}
In practice, we will truncate the expansions after, say, $p$ terms with the
value of $p$ to be determined later. We then define 
$\valpha \equiv (\alpha_{-p},\alpha_{-p+1},\dots,\alpha_0,\alpha_1,\dots,
\alpha_p)$ and 
$\vbeta \equiv (\beta_{-p},\beta_{-p+1},\dots,\beta_0,\beta_1,\dots,
\beta_p)$. 

\end{remark}	

\subsection{Multiple disks}
Suppose now that we have
$M$ well separated, identical dielectric disks randomly distributed in a 
homogeneous medium. Each disk is assumed to have radius $R$ and wavenumber $k_p$
and the background medium again has wavenumer $k_2$. For each individual particle, 
the analysis can be carried out as above. 
We
will denote by $\valpha^m$ the incoming coefficients and by $\vbeta^m$ 
the outgoing coefficients for the $m$-th particle. 
We have
\begin{equation}\label{scat}
\vbeta^m = S_p[\valpha^m], \mbox{ for } m = 1,\cdots,M.
\end{equation}
where $S_p$ denotes the truncated $(2p+1) \times (2p+1)$ 
scattering matrix acting on the truncated expansion.

The principle difference between the single particle and multi-particle scattering 
problem is that, in the latter case, the incoming field experienced
by each particle consists of two parts: the (applied) incident field $u^{inc}$ 
and the contribution to the scattered field $u^s$ from all of the other particles. 
In order to formulate the problem concisely, 
given the multipole expansion for disk $j$, 
we need some additional notation.

\begin{lemma} \cite{rok90}
Let disk $m$ be centered at $\mathbf{x_m}$
and let disk $l$ be centered at $\mathbf{x_l}$.
Then the multipole expansion 
\begin{equation}
\sum_{n=-\infty}^{\infty} \beta_n^m H_n(k_2 r_m)e^{i n\theta_m}
\end{equation}
induces a field on disk $l$ of the form
\begin{equation}
u = \sum_{n'=-\infty}^{\infty} \alpha^l_{n'} J_{n'}(k_2 r_l)e^{i n'\theta_l}
\end{equation}
where
\[
\alpha^l_{n'} = \sum_{n=-\infty}^{\infty} 
e^{-in(\theta_{lm}-\pi)} \beta^m_{n'-n}
				H_n(k_2 \| \mathbf{x_m} - \mathbf{x_l} \|).
\]
Here, $(r_m,\theta_m)$ and $(r_l,\theta_l)$ denote the polar coordinates
of a target point with respect to disk centers 
$\mathbf{x_m}$ and $\mathbf{x_l}$, respectively and
$\theta_{lm}$ denotes the angle between 
$(\mathbf{x_m} -\mathbf{x_l})$ and the $x$-axis.
\end{lemma}

\begin{remark}
We denote by 
$T^{jm}$ the translation operator that maps the outgoing coefficients $\vbeta^m$ 
from particle $m$ to the local expansion $\valpha^l$ centered at particle $l$. 
With this operator in place, the 
incoming coefficients $\valpha^m$ for the $m$-th particle is 
\begin{equation}\label{trans}
\valpha^m = \veca^m+\sum_{\substack{j = 1 \\ j\neq m}}^{M} T^{jm}\vbeta^j,
\end{equation}
where $\veca^m$ is the (truncated) local expansion \eqref{incident}
of the incident wave $u^{inc}$ on particle $m$.
$T^{jm}$ is referred to as the multipole-to-local (M2L) translation
operator \cite{rok90}.
\end{remark}

Combining eqs. \eqref{scat} and \eqref{trans}, one can easily eliminate the 
incoming coefficients $\valpha^m$ and obtain the following linear system that 
only involves the outgoing coefficients: 
\begin{equation}\label{scattmatrix}
\left(
\mathcal{S}^{-1}-\mathcal{T}\right)
\left[\begin{array}{c}
\vbeta^1 \\
\vbeta^2 \\
\vdots \\
\vbeta^M
\end{array} 
\right] 
=
\left[\begin{array}{c}
\veca^1 \\
\veca^2 \\
\vdots \\
\veca^M
\end{array} 
\right],
\end{equation}
where 
\begin{equation*}
\mathcal{S} = \left[\begin{array}{cccc}
S_p & & &\\
& S_p & &\\
& & \ddots &\\
& & & S_p
\end{array}
\right],
\quad
\mathcal{T} = \left[\begin{array}{cccc}
0 & T^{21} &\cdots &T^{M1}\\
T^{12} & 0 &\cdots &T^{M2}\\
\vdots&\vdots & \ddots &\vdots\\
T^{1M}&T^{2M} &\cdots & 0
\end{array}
\right].
\end{equation*}

The system 
\eqref{scattmatrix} can be solved iteratively, using GMRES \cite{GMRES1986}. 
Since each translation operator $T^{nm}$ is dense,
a naive matrix-vector product requires $O((M(2p+1))^2)$ operations, where $p$ is the 
order of the truncated expansion. 
FMM acceleration reduces the cost to $O(M(2p+1)^2)$ work,
for which we refer the reader to \cite{rok90,wideband}.
Further, \eqref{scattmatrix} has a simple diagonal preconditioner.
Multiplying through by the block diagonal matrix $\mathcal{S}$,
results in the preconditioned system matrix  $I-\mathcal{ST}$.
This significantly reduces the number of iterations. 

We now extend the multiple scattering approach to arbitrarily shaped particles.

\begin{remark}
It is worth emphasizing that the multiple scattering theory as discussed here 
is hardly new. We refer the reader to 
\cite{FL45,Gumerov2007206,Xu95} and the references therein. 
\end{remark}
   
\section{Wave scattering for arbitrarily shaped particles} 
\label{sec:arb}

When the dielectric inclusions are of arbitrary shape, multiple scattering
theory cannot be used quite so easily.  
Suppose, however, that an inclusion $\Omega$ is compactly supported with boundary
$\partial \Omega$ and that it is composed of a homogeneous material with wavenumber
$k_p$, as above.
Given the incident wave $u^{inc}$ and the boundary conditions 
\eqref{cont1}, \eqref{cont2}, the exterior scattered field $u^s$ and the 
field $u$ within $\Omega$ have the following representation \cite{Cot2}: 
\begin{align}
u^s&= \mathbf{S}^{k_2}\sigma+\mathbf{D}^{k_2}\mu, \mbox{ for } \mathbf{x}\in \Omega^c,\label{rep1}\\
u& = \mathbf{S}^{k_p}\sigma+\mathbf{D}^{k_p}\mu, \mbox{ for }\mathbf{x}\in \Omega, \label{rep2} 
\end{align} 
where $\mathbf{S}^{k}$ and $\mathbf{D}^{k}$ are the usual single layer and double layer potentials on $\partial \Omega$,
\begin{align}
\mathbf{S}^k \sigma & = \int_{\partial \Omega} G^{k}(\mathbf{x},\mathbf{y}) \sigma(\mathbf{y}) ds_{\mathbf{y}},\\
\mathbf{D}^k \mu &= \int_{\partial \Omega} \frac{\partial G^{k}(\mathbf{x},\mathbf{y})}{\partial n(\mathbf{y})}\mu(\mathbf{y})ds_{\mathbf{y}}.
\end{align}
$\sigma(y)$ and $\mu(y)$ are unknown charge and dipole densities that 
lie on the boundary $\partial \Omega$. 
We will need  the normal derivatives of $\mathbf{S}^k$ and $\mathbf{D}^k$ as well:
\begin{equation}
\mathbf{N}^k\sigma = \int_{\partial \Omega} \frac{\partial G^{k}(\mathbf{x},\mathbf{y})}{\partial n(\mathbf{x})}\sigma(\mathbf{y})ds_{\mathbf{y}}, \quad
\mathbf{T}^k \mu = \int_{\partial \Omega} \frac{\partial^2 G^{k}(\mathbf{x},\mathbf{y})}{\partial n(\mathbf{x}) \partial n(\mathbf{y})}\mu(\mathbf{y})ds_{\mathbf{y}}.
\end{equation}

By construction, the representations \eqref{rep1} and \eqref{rep2} satisfy 
the relevant Helmholtz equation in each domain. 
The single layer potential $\mathbf{S}^k$ is weakly singular and the value is 
well-defined for $\mathbf{x}\in \partial \Omega$. 
The operators $\mathbf{D}^k$ and $\mathbf{N}^k$ are define
on the boundary in the principal value sense (and have different limits 
when approaching the boundary from the interior and the exterior). 
The operator $\mathbf{T}^k$ is \textit{hypersingular} with its value on the 
boundary defined in the Hadamard finite part sense. 
For further details, we refer the reader to \cite{Cot2}. 

Enforcing the interface conditions \eqref{cont1}, \eqref{cont2} and taking 
appropriate limits \cite{Cot2} yields the following system of Fredholm
integral equations of the second kind:
\begin{align}
\mu+ [\mathbf{S}^{k_2}-\mathbf{S}^{k_p}]\sigma
+[\mathbf{D}^{k_2}-\mathbf{D}^{k_p}]\mu & = -u^{inc}, \label{intg1}\\
-\sigma+[\mathbf{N}^{k_2}-\mathbf{N}^{k_p}]\sigma
+[\mathbf{T}^{k_2}-\mathbf{T}^{k_p}]\mu & = -\frac{\partial u^{inc}}{\partial n}.
\label{intg2}
\end{align}

\begin{remark}
It is worth noting that, while $\mathbf{T}^k$ is hypersingular, the difference
kernel $\mathbf{T}^{k_2}-\mathbf{T}^{k_p}$ is only logarithmically singular
and compact as are all the other difference operators in \eqref{intg2}, at least
for smooth boundaries.
We use Nystr\"{o}m discretization for the system of equations based on the 
high order hybrid Gauss-trapezoidal rule of Alpert \cite{Alpert1999}.
In this paper, we restrict our attention to
smooth inclusions that are about one wavelength in size, so that
12 digits of accuracy are easily achieved with modest values of $N$ using the 
Gauss-trapezoidal rule for logarithmic singularities of order 16. 
We refer the reader to \cite{BG2011} and the references therein 
for further details.

The integral equation \eqref{intg2} was introduced in electromagnetics by M\"{u}ller
\cite{Muller}, and in the scalar case by Kress, Rokhlin, Haider, Shipman and 
Venakides \cite{Haider,KR78,rok2}.
\end{remark}

\subsection{The scattering matrix}

\begin{figure}[tbp]
\begin{center}
\includegraphics[width=80mm]{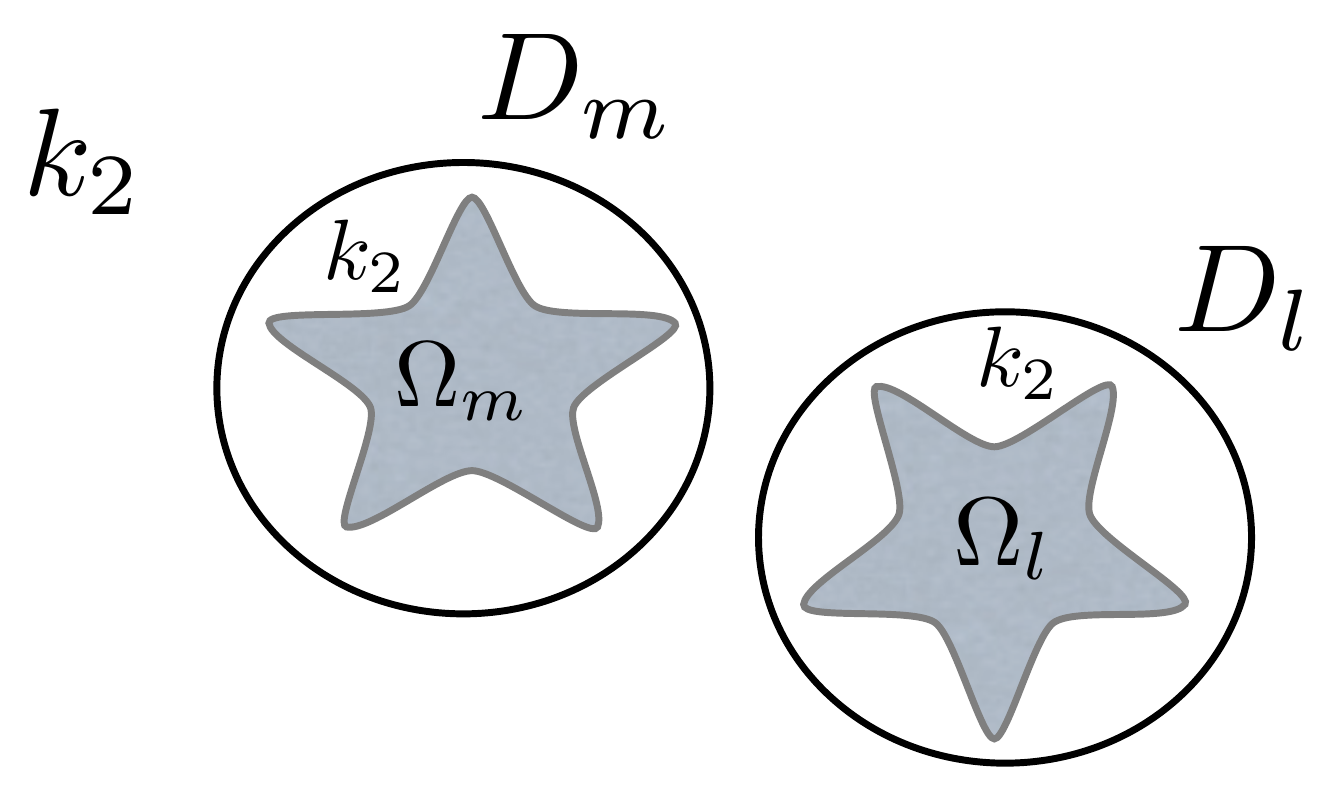}
\caption{Two inclusions and their enclosing disks. 
The scattering matrix $S_i$ for each inclusion $\Omega_i$ 
with wavenumber $k_p$ is defined as the map from an incoming
field on $D_i$ to the corresponding outgoing field. It is computed by solving
a sequence of boundary value problems on the inclusion itself in a precomputation
phase (see text).
In this paper, we assume that all the inclusions are identical but may be rotated,
as drawn here.
}
\label{scatdisk}
\end{center}
\end{figure}

Suppose now that we have $M$ inclusions $\Omega_1,\dots,\Omega_M$ that are 
identical up to rotation, and 
well separated in the sense that each inclusion $\Omega_i$ 
lies within a disk $D_i$ of radius $R$ so that the disks are not overlapping.
(see Fig. \ref{scatdisk}).

In that case, we can sample the incoming field on the disk $D_j$ rather than
$\Omega_j$ as 
\begin{equation}
u = \sum_{n=-p}^{p} \alpha_n J_n(k_2 r)e^{i n\theta} \, ,
\end{equation}
using a polar coordinate system centered on the disk $D_j$.

Let $\sigma_n$ and $\mu_n$ denote the solution to the 
integral equation \eqref{intg2} with right-hand side 
$u^{inc} = J_n(kr)e^{i n \theta}$,
$\frac{\partial u^{inc}}{\partial n} = kJ_n'(kr)e^{i n \theta}$.
We may then precompute the multipole expansion
from these source distributions 
\begin{equation}
u = \sum_{l=-p}^{p} \beta^n_l H_l(k_2 r)e^{i l\theta} \, ,
\end{equation}
where
\begin{equation}
\beta^n_l 
= \int_{\partial \Omega_j} 
[J_l(k_2|\mathbf{y}|) e^{- i l \theta_j(\mathbf{y})} \, 
\sigma_n(\mathbf{y})] + \, 
{\mathbf n} \cdot \nabla [J_l(w|\mathbf{y}|) e^{- i l \theta_j(\mathbf{y})}
\mu_n(\mathbf{y})] \, 
ds_{\mathbf{y}} \, .
\end{equation}
Here, $\mathbf{y}$ is the location of a point on 
$\partial \Omega_j$ with respect to the center of disk $D_j$ and
$\theta_j(\mathbf{y})$ is the polar angle subtended with respect to
the center of disk $D_j$.
The formula for $\beta_l$ is standard 
\cite{rok90,wideband} and derived from the Graf addition theorem
\cite{Hand2010}. 
\begin{definition}
As before,
the mapping between the incoming coefficients $\{\alpha_n\}$ and 
outgoing coefficients $\{\beta_n\}$ is referred as the scattering matrix for 
the inclusion $\Omega_j$ and denoted by $S_j$.
\end{definition}

The reason for permitting a different scattering matrix for each inclusion
is that the $\Omega_j$ may be distinct in terms of geometry or
dielectric properties. For the sake of simplicity, we assume here that
the wavenumbers are the 
same in each inclusion and that the shapes are the same up to rotation.
This permits us to solve only $2p+1$ integral equations on a single
prototype inclusion in the enclosing disk. 
The scattering matrix for each rotated copy is then trivial
to construct.
Moreover, we can easily store the 
densities $\sigma_n$ and $\mu_n$, since this requires only 
$O(2N(2p+1))$ storage, where $N$ is the number of points used to discretize 
the boundary $\partial \Omega$. The amount of memory required to store the
scattering matrix is $O((2p+1)^2)$. For modest values of $N$, as is the case
in the present paper,
we compute the $LU$ factors of the integral equation system matrix 
corresponding to \eqref{intg1}, \eqref{intg2} only once,
at a cost of $O(N^3)$ work. Each right-hand side 
corresponding to 
$u^{inc} = J_n(k_2R)e^{i n \theta}$ and 
$\frac{\partial u^{inc}}{\partial n} = k_2 J_n'(k_2 R)e^{i n \theta}$ 
can then be solved for $n = -p,\dots,p$ at a total cost of 
$O(N^2(2p+1))$ work.

\subsection{Multiple scattering}

If we were interested in solving the multiple scattering problem in an 
infinite medium, we could now proceed as in the previous section.
The number of degrees of freedom is only $2p+1$ per inclusion rather than
$N$ points per inclusions (the number needed to discretize the domain boundaries
$\partial \Omega_j$).
For complicated inclusions, this permits a vast reduction
in the number of degrees of freedom required and forms the basis for 
the FMPS method \cite{GG2013}. Moreover, the block-diagonal preconditioned
multiple scattering equations are much better conditioned than the 
integral equation \eqref{intg1}, \eqref{intg2} itself and FMM acceleration
is particularly fast in this setting.

\begin{remark}
Extending the method to more than one type of substructure is straightforward as 
long as the assumption that the enclosed circles are well separated still holds. 
The additional cost is the bookkeeping for different scattering matrices of these 
substructures.
\end{remark}
 
\section{Multi-particle scattering in a layered medium}

To this point, we have discussed the layered medium and multiple scattering
problem spearately. For the full problem, we now assume that
multiple inclusions have been placed in the middle of a three-layered medium. 
We assume that the inclusions are well separated, so that the multiple
scattering formalism applies within the layer.
Then, we may write 
\begin{align}
u_1(\mathbf{x}) &= G_{k_1}(\mathbf{x},\mathbf{x}_0) + u_1^s  \nonumber \\
u_2(\mathbf{x}) &= u_2^t + u_2^b + \sum_{j=1}^M \sum_{n=-p}^p 
\beta_n^m H_n(k_2r_m) e^{i n \theta_m} \label{fullrep} \\
u_3(\mathbf{x}) &= u_3^s \nonumber 
\end{align}
where $u_1$ and $u_3$ denote the fields in the top and bottom half spaces and
$u_2$ denotes the field in the central layer exterior to the scattering disks $D_j$.
$u_1^s$, $u_2^t$, $u_2^b$, and $u_3^s$ are the Sommerfeld integrals from 
Section \ref{sec2}. Once $u_2$ is known, the field within the scattering disks
and the inclusions themselves is easily obtained.

It remains to discuss the discretization of the Sommerfeld integral,
and the setup of the global linear system for the unknowns
$\sigma_1$, $\sigma_2^+$, $\sigma_2^-$, $\sigma_3$, and 
$\{ \vbeta_m, m = 1\dots,M \}$.

\subsection{Evaluation of the Sommerfeld integral}

\begin{figure}[tbp]
\begin{center}
\includegraphics[width=80mm]{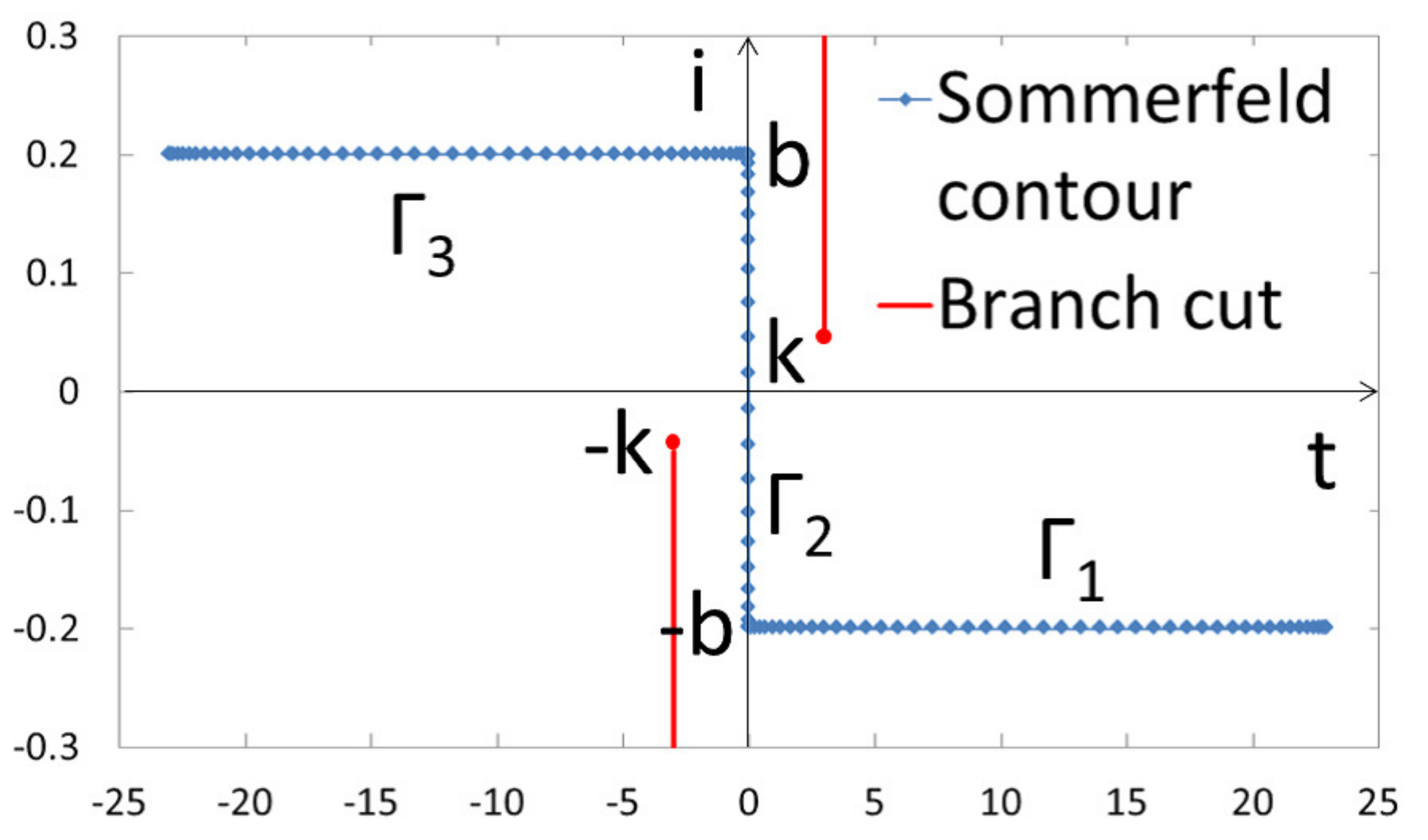}
\caption{The Sommerfeld contour in the complex $\lambda$ plane: Each segment in the 
contour is discretized using Gauss-Legendre quadrature. 
The branch cut (shown in red) points upward from $k$ and downward from $-k$.}
\label{fig:Sommerfeldcontour}
\end{center}
\end{figure}

Let us consider the function $u_2^t$ defined 
by \eqref{layer22}.
Its computation is a standard problem in acoustic and 
electromagnetic scattering and often handled by contour deformation.
It is typical to deform the integration contour by pushing it from the real line
into the second and fourth quadrants of the complex $\lambda$-plane in order to 
avoid the square root singularities in the integrand. 
One option is to use a 
hyperbolic tangent contour \cite{BG2011}, which yields spectral accuracy with
the trapezoidal rule and is extremely efficient.  
In our numerical simulation, we have chosen to use the piecewise smooth contour 
shown in Fig. \ref{fig:Sommerfeldcontour} instead. This is slightly less efficient,
but will permit us to evaluate the Sommerfeld integral using the non-uniform FFT,
as explained further below.
The contour consists of three segments: $\Gamma_1$, $\Gamma_2$ and $\Gamma_3$, where
\begin{equation}\label{sommerseg}
\left\{\begin{array}{ccc}
\Gamma_1: &t - ib, & t\in (0,\infty), \\
\Gamma_2: & it,  & t \in [-b,b], \\
\Gamma_3: &t + ib, & t\in (-\infty,0).
\end{array}\right.
\end{equation}
The branch cuts for the square root in the integrand are chosen to ensure
that waves are decaying away from the interface.
(up at $k$ and down at $-k$ as shown in Fig. \ref{fig:Sommerfeldcontour}).

We truncate 
$\Gamma_1$ and $\Gamma_3$ at a point $t_{max}>0$, where the 
integrand of $u_2^t$ has decayed to a user-specified tolerance.
Fortunately, the decay in the integrand is exponential 
once $\lambda$ exceeds $k_2$.
(The precise rate of decay depends on the distance from the interface
of the scattering disks and the point source generating the incoming field.)
We let $N_S$ denote the number of points used in the quadrature for the 
Sommerfeld contour and note that
each discretization point $\lambda_j$ on the contour corresponds to a plane wave. 
We use the same contour and the same $N_S$ values $\{ \lambda_j \}$ for each
of $u_1^s$, $u_2^t$, $u_2^b$, and $u_3^s$.

\subsection{The full linear system}

Let us denote by $\vsigma$ the discretized densities on the dielectric layers, 
$\vsigma = [\vsigma_1,\vsigma_2^+,\vsigma_2^-,\vsigma_3]^T$, and by $\vbeta$ 
the multipole coefficients for all $M$ particles in the central layer. 
Each of $\vsigma_1$, $\vsigma_2^+$, $\vsigma_2^-$,  and 
$\vsigma_3$ is of length $N_S$ and the 
full linear system for multiple scattering in the layered medium takes the form of
a block $2\times 2$ linear system:
\begin{equation}\label{finalsys}
\left[\begin{array}{cc}
A & B \\
C & D 
\end{array}\right]
\left[\begin{array}{c}
 \vsigma \\
 \vbeta
\end{array} \right]=
\left[\begin{array}{c}
b \\
0
\end{array}\right].
\end{equation}  

$A$ itself is block diagonal $4N_S \times 4N_S$ matrix with $4 \times 4$ blocks 
of the form $A_\lambda$ in \eqref{sommerfeldsys}, each such block corresponding
to a distinct $\lambda_j$ in the contour integral discretization.
The right-hand side component $b$ is simply the right-hand side of
\eqref{sommerfeldsys} for each such $\lambda_j$. 
The matrix $D =\mathcal{ST}-\mathcal{I}$ is simply the
multiple scattering system for the particles from 
\eqref{scattmatrix}.
The off-diagonal blocks $B$ and $C$ are more complicated.
$B$ is a matrix that translates the multipole expansion coefficients
to a Sommerfeld representation on the upper and lower interfaces of
the layered medium,  while
$C$ requires the evaluation of the Sommerfield integral contributions from
the interfaces in terms of incoming local expansions on the scattering disks
themselves. 
We turn now to the efficient application of the matrices $B$ and $C$.

\subsubsection{The Sommerfeld-to-local operator}
\label{sec5}
A straightforward mechanism to map from the $\vsigma$ variables to 
local expansions on the $M$ disks is to use the 
Jacobi-Anger formula \cite{Hand2010}.
\begin{lemma}
\label{lem11}
Given $r\in \mathbb{R}$, $k\in\mathbb{C}$, we have
\begin{equation}
e^{ikr\cos{\theta}} = \sum_{n=-\infty}^{\infty} i^nJ_n(kr)e^{in\theta}.
\end{equation} 
\end{lemma}
Suppose now that we want to compute the contribution from $\sigma_2^+$ to a 
local expansion on a disk centered at $(x_1,y_1)$.
Using Lemma \ref{lem11}, it is easy to see that 
\begin{equation}\label{JA2}
e^{\sqrt{\lambda_j^2-{k_2}^2}y+i\lambda_j(x-x_0)} = e^{\sqrt{\lambda_j^2-{k_2}^2}y_1+i\lambda_j(x_1-x_0)} \sum_{n=-\infty}^{\infty} i^nJ_n(k_2r)e^{in(\phi+\theta)},
\end{equation}
where $\phi = \arccos(\lambda_j/k_2)$, 
$\theta = \arccos((x-x_1)/r)$ and $r = \sqrt{(x-x_1)^2+(y-y_1)^2}$.
The analogous formula can be obtained for 
the contribution from $\sigma_2^-$.

The cost of using 
formula \eqref{JA2} to compute the action of the $C$ block in the system matrix above
is clearly $O(MN_S(2p+1))$, where $M$ 
denotes the number of particles and $N_S$ the number of discretization points $\lambda_j$ in
the Sommmerfeld contour and $p$ is the order of the expansions used in the multiple scattering
representation. This is quite acceptable when either $N_S$ or $M$ is small. 
For high frequency problems with many inclusions, where $k_2$ is large and $N_S = O(k_2)$, 
we have developed a more efficient scheme, based on the nonuniform FFT (NUFFT).

\subsubsection{The Sommerfeld-to-local operator using the NUFFT}
\label{sec53}

Instead of mapping the contribution from the Sommerfeld integral to each disk separately,
we seek a fast algorithm for evaluating the integral on a grid of points in the central layer,
after which we can use high order interpolation to get the desired local expansion.

Restricting our attention to $u_2^t$ for a fixed value of $y$, we have 
\begin{equation}\label{fft1}
 \frac{1}{4\pi}\int_{\Gamma_1}\frac{e^{\sqrt{\lambda^2-{k_2}^2}y}}{\sqrt{\lambda^2-{k_2}^2}} e^{i\lambda(x-x_0)}\sigma_2^+(\lambda)d\lambda = 
\frac{1}{4\pi}e^{b(x-x_0)} \int_0^{t_{max}}g(t)e^{itx}dt, 
\end{equation} 
where
\[
g(t) = \frac{e^{\sqrt{(t-ib)^2-k_2^2}y}}{\sqrt{(t-ib)^2-k_2^2}}
\, e^{-it x_0} \, \sigma_2^+(t-ib).
\]

Note now that the integral on the right-hand side of \eqref{fft1} is a finite Fourier transform.
If we could compute it rapidly, we would have an efficient method for evaluating the Sommerfeld
integral at a fine grid in the $x$ variable for a fixed $y$.
The discretization points in $t$, however, lie at Gauss-Legendre nodes, so the FFT itself does
not apply. Fortunately, the nonuniform FFT (NUFFT) of Dutt and Rokhlin \cite{DR93,Dutt199585}
permits this to be done in nearly linear time.
In our numerical simulations, we use the version discussed in \cite{GL2004,Lee20051}. 
The analogous method permits the rapid evaluation of the 
Sommerfeld integral on the contour $\Gamma_3$. 
For the integral on $\Gamma_2$, the NUFFT cannot be applied, but only a few 
discretization points are required, so we evaluate that contribution directly.

To provide rapid access to the field induced by the Sommerfeld integral at any location in the
central layer, we superimpose on it a grid of $n_1 \times n_2$ boxes that
contain all of the $M$ scattering disks. In each such box, we construct a tensor product
$m_1 \times m_2$ Chebyshev mesh, which will permit $q$th order local interpolation by 
barycentric interpolation \cite{BT2004}.
The cost for evaluation at all grid points is $O\left( (n_2 m_2) \,  (n_1 m_1 + N_S) \log(n_1m_1  + N_S)
\right)$ operations, using the NUFFT for each of the distinct $n_2 m_2$ locations in $y$.

Consider now one of the scattering disks $D_j$ of radius $R$. 
If we discretize the boundary of the disk
using $2p+1$ equispaced points, evaluation of the induced field at each of the points 
requires $O(m_1m_2)$ operations, for a net cost of $O(m_1m_2(2p+1))$ work.
An FFT of order $(2p+1)$ converts these values into their Fourier transforms, which we denote by
$\{u_n\}$, for $n = -p,\cdots, p$. From this, the $n$-th term $a_n$ in the 
incoming $J$-expansion is simply
\begin{equation}\label{proj}
a_n = \frac{u_n}{J_n(k_2 R)}.
\end{equation}

\begin{remark}
The formula \eqref{proj} will fail if the value $k_2R$ is a zero of the function $J_n$ for any 
$n$ from $-p\,\dots,p$. This can be avoided if we also compute the 
normal derivative of the Sommerfeld integral on the boundary of each scattering disk.
If we denote by $\{u'_n\}$ the Fourier coefficient of the normal derivative, it is
easy to see that
\begin{equation}
a_n = \frac{u_nJ_n(k_2R)+u'_nkJ'_n(k_2R)}{J^2_n(k_2R)+(kJ'_n(k_2R))^2}, \mbox{ for }j = -p,\cdots,p.
\end{equation}
The evaluation of the gradient of the Sommerfeld integral can be computed by an obvious modification 
of the formula \eqref{fft1} or (with a reduction in order) 
by computing the gradient of the tensor product Chebyshev series discussed above.
\end{remark}

In summary, it requires
$O(Mm_1m_2(2p+1))$ operations to interpolate the field values on each of the $M$ scattering
disks and $O(M(2p+1)\log(2p+1))$ operations to obtain the coefficients of the $J$-expansions. 
This completes the computation of the $C$ block in the system matrix.

\subsection{The multipole-to-Sommerfeld operator}

The off-diagonal $B$ block in \eqref{finalsys} requires a formula for recasting the 
multipole expansion to the corresponding Sommerfeld representation on either the upper or lower
interface of the layered medium.
More precisely, each $H$-expansion in the central layer, centered on disk $D_j$ with center
$(x_j,y_j)$ has
a spectral representation on the upper layer $y= 0$ and the lower layer $y= -d$ of the form:
\begin{align}
u_{j}^t &= \frac{1}{4\pi}\int_{-\infty}^{\infty}\frac{1}{\sqrt{\lambda^2-{k_2}^2}} 
e^{i\lambda(x-x_0)}\sigma_{mp}^{+}(\lambda)d\lambda, \label{layer211}\\
u_{j}^b &= \frac{1}{4\pi}\int_{-\infty}^{\infty}\frac{1}{\sqrt{\lambda^2-{k_2}^2}} 
e^{i\lambda(x-x_0)}\sigma_{mp}^{-}(\lambda)d\lambda, \label{layer221}
\end{align}
respectively.

The formulae for $\sigma_{mp}^{+}(\lambda)$ and $\sigma_{mp}^{-}(\lambda)$ follow directly from
the following theorem.

\begin{theorem}\label{thm11} \cite{Cheng2006616}
Let $(x_j,y_j)$ denote the center of a multipole expansion in the central layer, 
with $-d < y_j <0$ and let $(r,\theta)$ denote the polar coordinates of a target point 
with respect to that center. Then,
on the upper interface,
\begin{equation}\label{planeexp2}
H_n(k r)e^{i n \theta} = \frac{(-1)^n}{4\pi}\int_{-\infty}^{\infty} 
\frac{e^{\sqrt{\lambda^2-k^2}y_j}}{\sqrt{\lambda^2-k^2}}e^{i\lambda (x-x_j)}  
\bigg( \frac{\sqrt{\lambda^2-k^2}+k^2}{k^2}\bigg)^n d\lambda,
\end{equation}
and on the lower interface,
\begin{equation}\label{planeexp1}
H_n(k r)e^{i n \theta} = \frac{(-1)^n}{4\pi}\int_{-\infty}^{\infty} 
\frac{e^{-\sqrt{\lambda^2-k^2}(d+y_j)}}{\sqrt{\lambda^2-k^2}}e^{i\lambda (x-x_j)} 
\bigg( \frac{\sqrt{\lambda^2-k^2}-k^2}{k^2}\bigg)^n d\lambda.
\end{equation}
\end{theorem}

Each multipole coefficient in the expansion about disk $D_j$ contributes to each of the 
$N_S$ discretization points in the Sommerfeld integrals, requiring a total of 
$O\left( (2p+1)N_SM \right)$ work. This, then, is the cost of applying the $B$ block 
of the system matrix directly.

\subsubsection{The multipole-to-Sommerfeld operator using the NUFFT}

Because of the computational complexity of applying the $B$ block in the manner described above,
it is important to develop a fast algorithm for the case where $M$ and $N_S$ are large.
We do so by essentially inverting the method of section \ref{sec53}.
Assume first that all the centers of the $H$-expansions lie at the nodes of a uniform grid in the central
layer and let us consider the contributions from the $n$th mode at each such grid point $(x_l,y_j)$ 
for a fixed horizontal line $y=y_j$.
If there are $n_1$ such expansion centers, with $x$ coordinates $\{x_l\}$, $l=1,\cdots,n_1$,
and we denote by 
${a_n^l}$ the coefficient for the $n$th mode of the $H$-expansion at location $(x_l,y_j)$,
then the total contribution to the induced spectral coefficient
$\sigma_{mp}^+(\lambda_j)$ on the top layer is given by 

\begin{eqnarray}\label{invfft}
\{\sigma_{mp}^+(\lambda_j)\}_n &:=& e^{\sqrt{\lambda_j^2-k_2^2}y_j}
\bigg( \frac{\sqrt{\lambda_j^2-k^2}+k^2}{k^2}\bigg)^n 
 \sum_{l=1}^{n_1} a_n^l e^{-i\lambda_j(x_j-x_0)} \\
\{\sigma_{mp}^-(\lambda_j)\}_n &:=& e^{-\sqrt{\lambda_j^2-k_2^2}(d+y_j)}
\bigg( \frac{\sqrt{\lambda_j^2-k^2}-k^2}{k^2}\bigg)^n 
 \sum_{l=1}^{n_1} a_n^l e^{-i\lambda_j(x_j-x_0)} \nonumber
\end{eqnarray}
 
The formulae \eqref{invfft} imply that for each row, one can use the NUFFT to compute
the induced coefficients for each discrete quadrature node $\lambda_j$ on $\Gamma_1$ or $\Gamma_3$. 
As above, we use direct computation for the contributions to discretization nodes on $\Gamma_2$. 
In the general case, the 
centers of the $H$-expansions are not aligned on a grid, but we can first shift the center
of each $H$-expansion to the nearest grid point, using the multipole-to-multipole
translation operator \cite{rok90,wideband} based on the Graf addition theorem \cite{Hand2010}. 
After $M$ such shifts, we may apply the transformation of \eqref{invfft}.

The total computational cost is $O(M (2p+1)^2)$ for shifting all the $H$-expansions and 
 $O\left( n_2 \, (2p+1) \, (n_1 + N_S) \log(n_1 + N_S) \right)$ for the  NUFFT-based work
(see Table \ref{tab1}). The merits of the NUFFT-based schemes would become more apparent for larger $N_S$.

\begin{table}[htbp]
\begin{center}
\caption{Comparison of CPU time in seconds for the Sommerfeld-to-local and multipole-to-Sommerfeld
operators, unsing both the direct and NUFFT-based schemes (see text).
The Sommerfeld contour is discretized with $500$ Gauss-Legendre points
($240$ points for $\Gamma_1$ and $\Gamma_3$, with $20$ points for $\Gamma_2$).} 
\label{tab1}
\begin{minipage}{\linewidth}
\begin{center}
\subcaption{Computation of the Sommerfeld-to-local operator}
\begin{tabular}{c|cccc} \hline
Number of scatterers & 100 & 500 & 1,000 & 5,000\\ \hline
Direct method & 2.21e-2 & 9.90e-2 & 2.06e-1 & 9.86e-1\\ \hline
NUFFT & 1.31e-1 & 1.77e-1 & 2.28e-1 & 4.93e-1\\ \hline
\end{tabular}
\end{center}
\vspace{0.3cm}
\end{minipage}
\begin{minipage}{\linewidth}
\begin{center}
\subcaption{Computation of the multipole-to-Sommerfeld operator}
\begin{tabular}{c|cccc} \hline
Number of scatterers & 100 & 500 & 1,000 & 5,000\\ \hline
Direct method & 3.49e-2 & 1.82e-1 & 3.60e-1 & 1.80\\ \hline
NUFFT & 6.05e-2 & 1.38e-1 & 1.61e-1 & 2.77e-1\\ \hline
\end{tabular}
\end{center}
\end{minipage}
\end{center}
\end{table}

\subsection{Iterative solution of the system matrix}

We will solve equation \eqref{finalsys} using the iterative method GMRES \cite{GMRES1986}. 
However, the unknowns $\vsigma$ and $\vbeta$ may be poorly scaled with respect to each other.
However, $A$ is block diagonal, as noted above, with simple $4 \times 4$ blocks. 
Thus, we first invert $A$ directly and use GMRES on the Schur complement of \eqref{finalsys}. 
In other words, we solve the system
\begin{equation}\label{schur}
[D-CA^{-1}B][\vbeta] = -CA^{-1}b
\end{equation}
instead.
This is much better conditioned and involves only the $\vbeta$ unknowns.
The Schur complement formalism has a simple physical interpretation:
it is, in essence, a reformulation of the scattering problem using the layered medium 
Green's function.
  
\section{Numerical experiments}

In this section, we illustrate the performance of our algorithm with three 
examples. For simplicity, we use a single class of inclusions, 
parametrized by 
\begin{equation}\label{par1}
\left\{\begin{array}{rcl}
x &=& (a_1+a_2\cos(a_3 t))\cos(t),\\
y &=& (a_1+a_2\cos(a_3 t))\sin(t),
\end{array}\right.\mbox{ for } 0\leq t < 2\pi.
\end{equation}
As discussed in section \ref{sec:arb},
inclusions with more complicated boundaries do not introduce any essential 
difficulty in our scheme except that the precomputation of the scattering matrix 
is a little more involved, particulalry if corners are present 
\cite{BRS2010,Helsing2008}.

Given a fixed $a_1$, $a_2$ and $a_3$, multiple copies of the inclusion are randomly 
distributed in the central layer of the medium with random orientations. 
To ensure the inclusions are well separated but confined in a fixed region, 
we use a \textit{bin sorting} algorithm to construct the random distribution. 
We begin with inclusions located on a regular grid and then perturb their positions randomly,
accepting the random move if the inclusion remains inside the region and 
stays well separated from the others. Several such sweeps are carried out to randomize the
positions further.

We have not, as yet, specified the parameter $N_S$ used to discretize the Sommerfeld integral 
in \eqref{sommerseg}. While special techniques have been developed by many authors to 
handle sources near the interface (see \cite{BG2011,ChoCai,ONeil2014}), we simply
assume that the source defining the incoming field is at least 0.2 wavelengths from the top interface.
More precisely, in our examples, the source point in the top layer is placed
at $(1,1)$ (which is roughly $0.2$ wavelengths away for wavenumber $k_1=1$). 
We also assume that the nearest the inclusions get to either one of the interfaces in 
the layered medium is at least $0.5$ wavelengths. 
Under these assumptions, we let $t_{max} = \max\{|k_1|,|k_2|,|k_3|\}+20$, $b = 0.2$, 
and discretize $\Gamma_1$ and $\Gamma_3$ using $240$ Gauss-Legendre points and $\Gamma_2$ by $20$ 
Gauss-Legendre points. This is sufficient to achieve about $10$ digits of accuracy.

All computations are carried out using a 2.3GHz Intel Core i5 laptop, with 4GB RAM.

\subsection{Example 1: scattering from large numbers of inclusions }
\label{NumericalExperiment1}

\begin{figure}[tbp]
\begin{center}
\includegraphics[width=140mm]{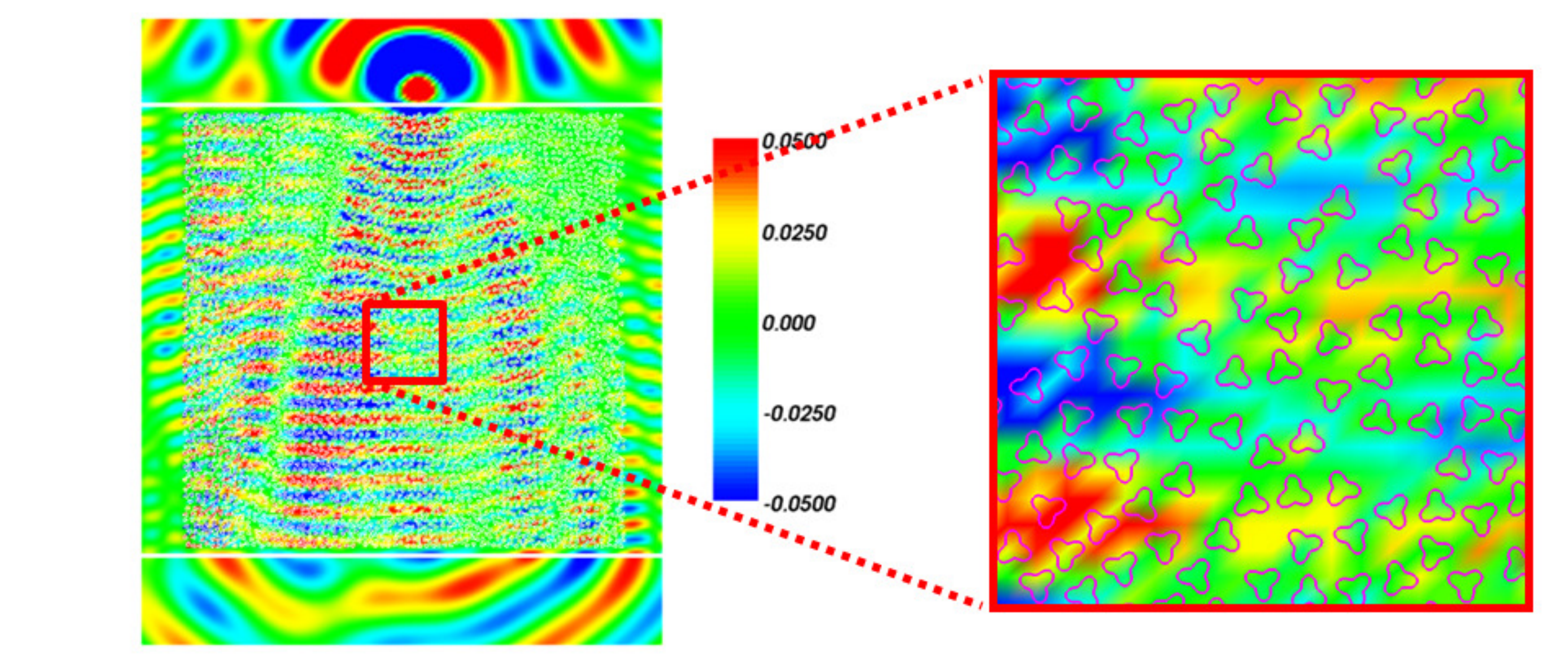}
\caption{Real part of the total field with $5,000$ dielectric inclusions randomly distributed 
in a three-layered medium. The wavenumber for each particle is $k_p = 2.0$ and the 
wavenumbers for the three layers are $k_1=1.0$, $k_2=3.0$, $k_3=1.0$. 
The size of each particle is approximately $0.1$ wavelength for the wavenumber $k_2$.
}
\label{fig:NumericalResult1}
\end{center}
%
\begin{minipage}{0.5\linewidth}
\begin{center}
\includegraphics[width=70mm]{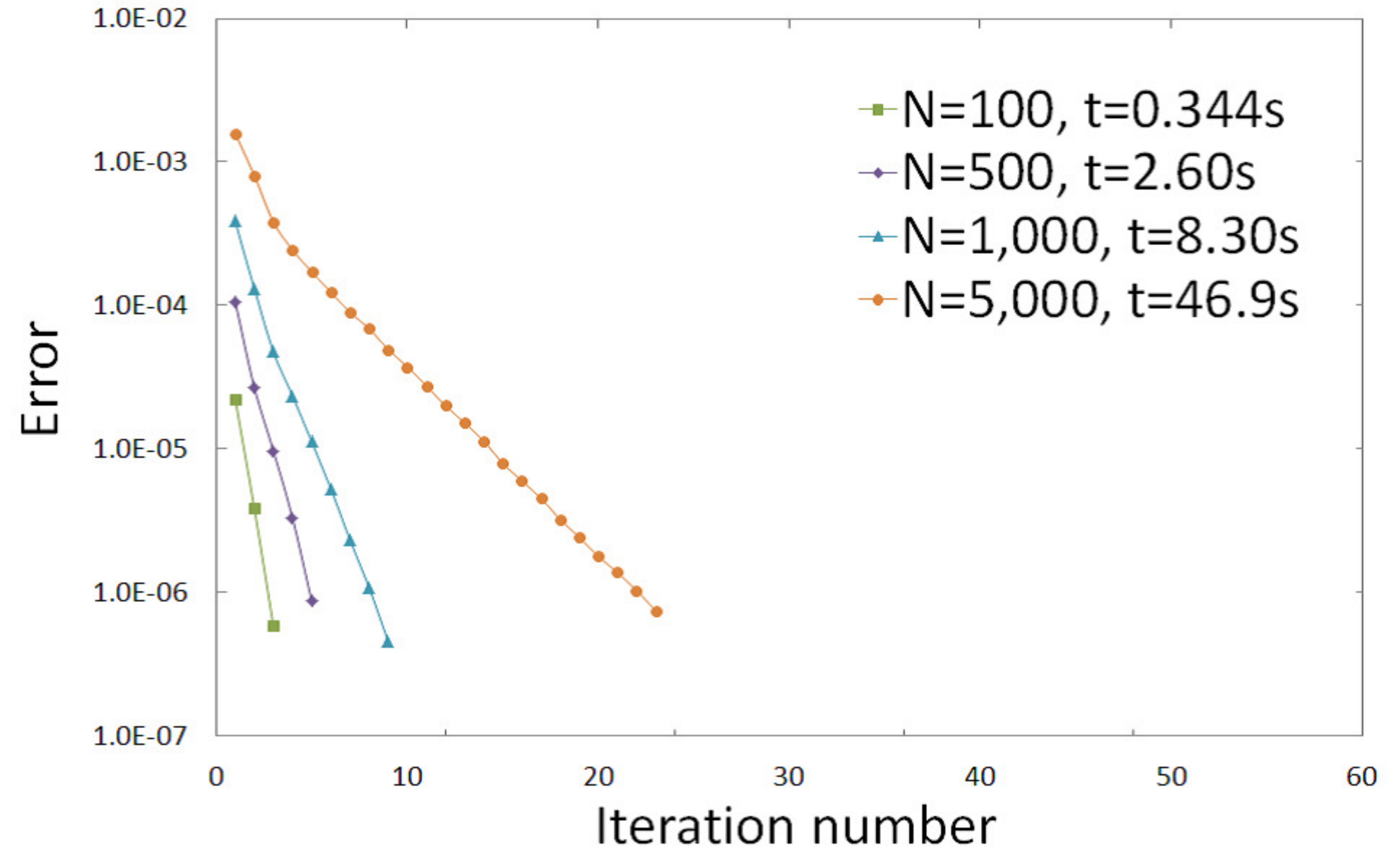}
\subcaption{}
\end{center}
\end{minipage}
\begin{minipage}{0.5\linewidth}
\begin{center}
\includegraphics[width=70mm]{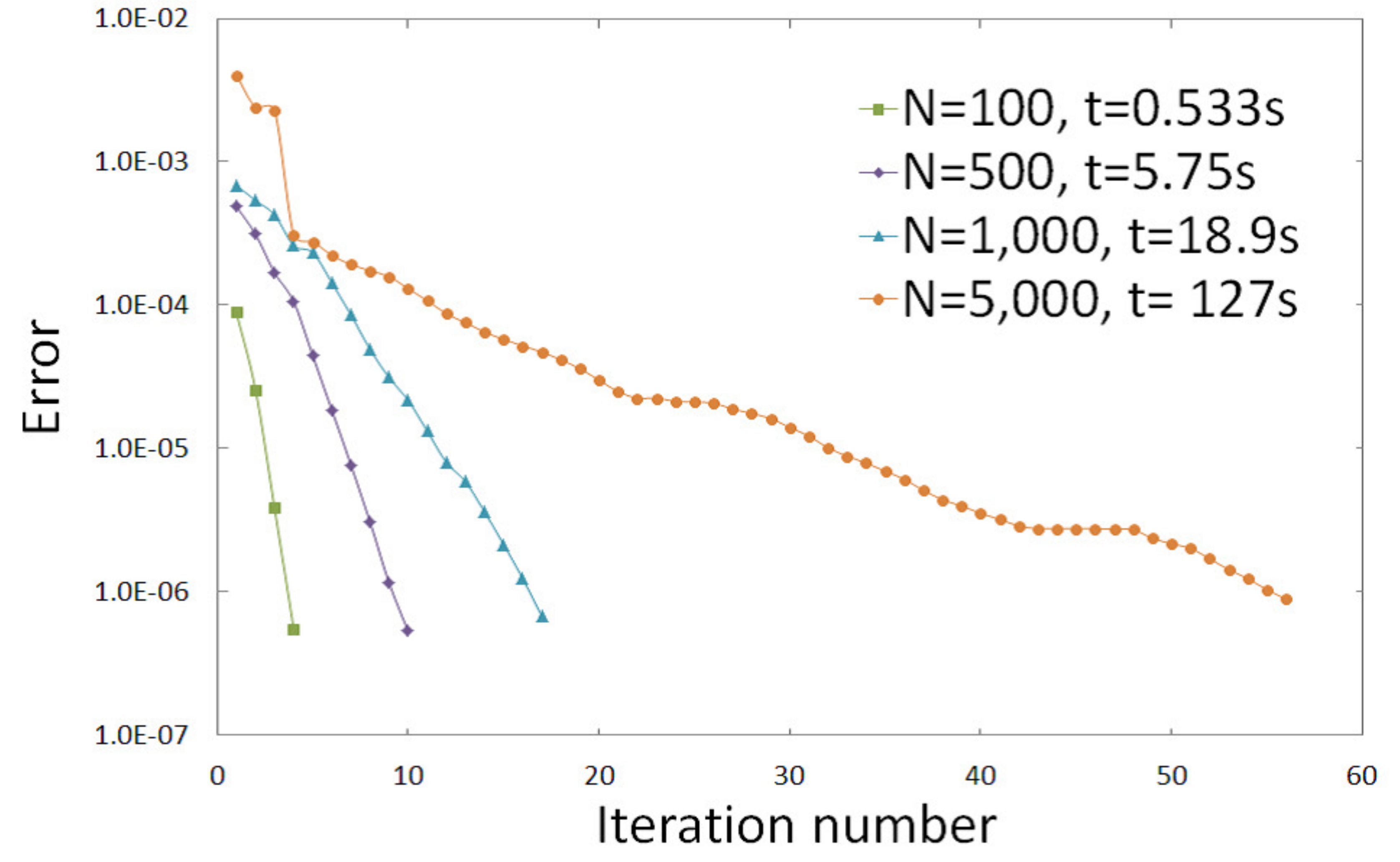}
\subcaption{}
\end{center}
\end{minipage}
\begin{center}
\caption{Convergence behavior of GMRES and the CPU time required
for various numbers of inclusions embedded in either (a) free space or 
(b)a three layered medium. For (a), we set $k_1=k_2=k_3=3.0$ and for (b),
we set $k_1=1.0$, $k_2=3.0$, $k_3=1.0$.}
\label{fig:NumericalResult1-1}
\end{center}
\end{figure}

In our first example, we consider the scattering of inclusions defined by parameters
$a_1=0.12,\ a_2=0.04$, and $a_3=3$ in eq. \eqref{par1} with wavenumber $k_p = 2.0$. 
To obtain the scattering matrix with $p = 10$, we solve the integral equation 
\eqref{intg1} and \eqref{intg2} by discretizing the boundary of the particle using
$N=300$ equispaced points.  We assume the wavenumbers of the layered medium are given by 
$k_1=1.0,\ k_2=3.0,\ k_3=1.0$. The thickness of the central layer is determined by the 
parameter $d=32$. We consider distributions of $M= 100,\ 500,\ 1,000,\ 5,000$  inclusions
and solve the mulitple scattering problem using GMRES with FMM acceleration. 
We terminate the iteration once the residual is less than $10^{-6}$. 
Results are presented in Fig. \ref{fig:NumericalResult1} and \ref{fig:NumericalResult1-1}. 

Fig. \ref{fig:NumericalResult1} shows the total field in the case $M= 5,000$. 
The field distortion due to the inclusions is apparent. 
It requires $127$s to achieve 6 digits of accuracy. 
Fig. \ref{fig:NumericalResult1-1} shows the convergence behavior of GMRES as the number of inclusions
is increased as well as the total CPU time. Clearly, more iterations are required for 
larger numbers of particles. Nevertheless, the time scales roughly linearly with the number
of particles. In Fig. \ref{fig:NumericalResult1-1}(a), we study the convergence rate 
when the background is homogeneous, by setting the material parameters to be the same for the three 
layers ($k_1 = k_2 = k_3$). As expected, convergence is more rapid than when the inclusions
are embedded in a true layered medium, because of the multiple reflections from the 
interfaces themselves.

\subsection{Example 2: scattering in high contrast materials}
 \label{NumericalExperiment2}

\begin{figure}[tbp]
\begin{center}
\includegraphics[width=140mm]{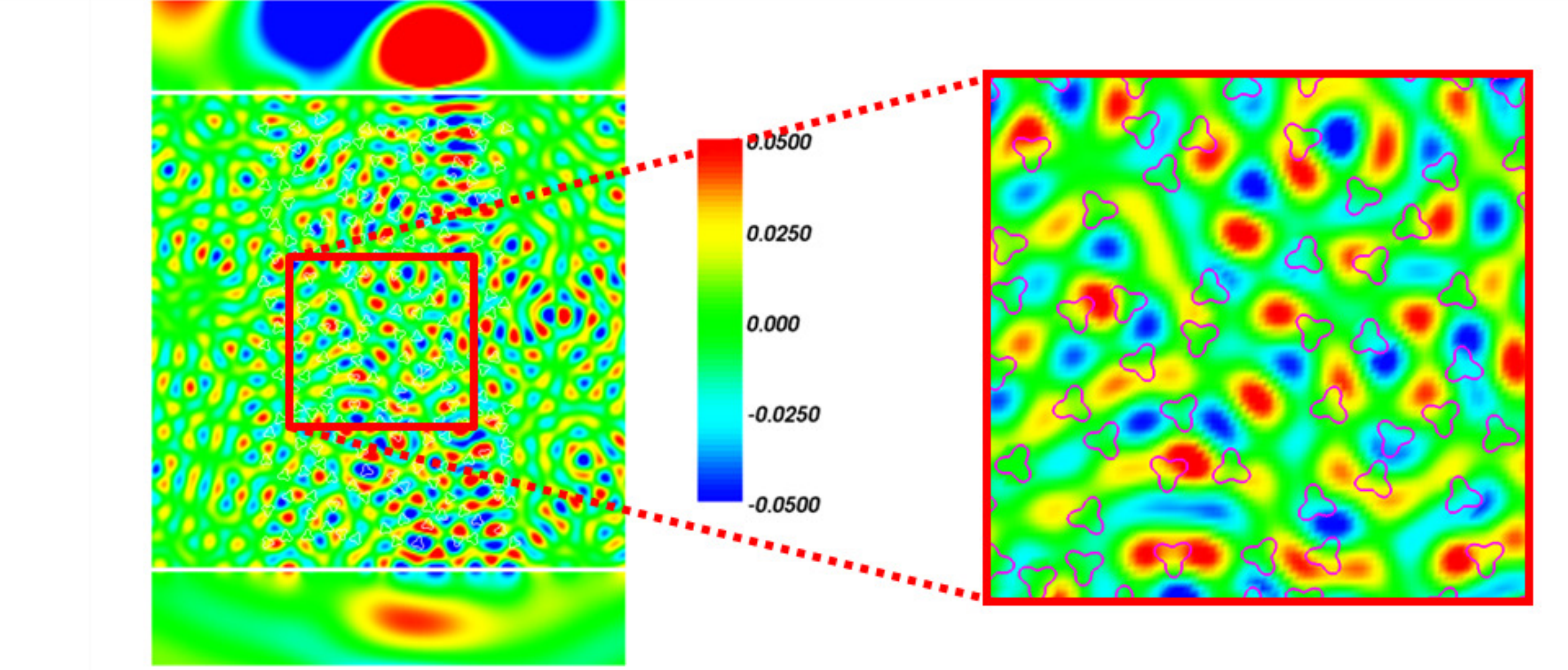}
\caption{Real part of the total field for $200$ dielectric inclusions distributed in 
a three layer medium with wavenumbers $k_1=1.0$, $k_2=10.0$, $k_3=1.0$.
For each inclusion, the wavenumber is $k_p = 2.0$. 
The inclusions are approximately $0.3$ wavelength in size for the 
wavenumber $k_2$.}
\label{fig:NumericalResult2}
\end{center}
%
\begin{minipage}{0.5\linewidth}
\begin{center}
\includegraphics[width=70mm]{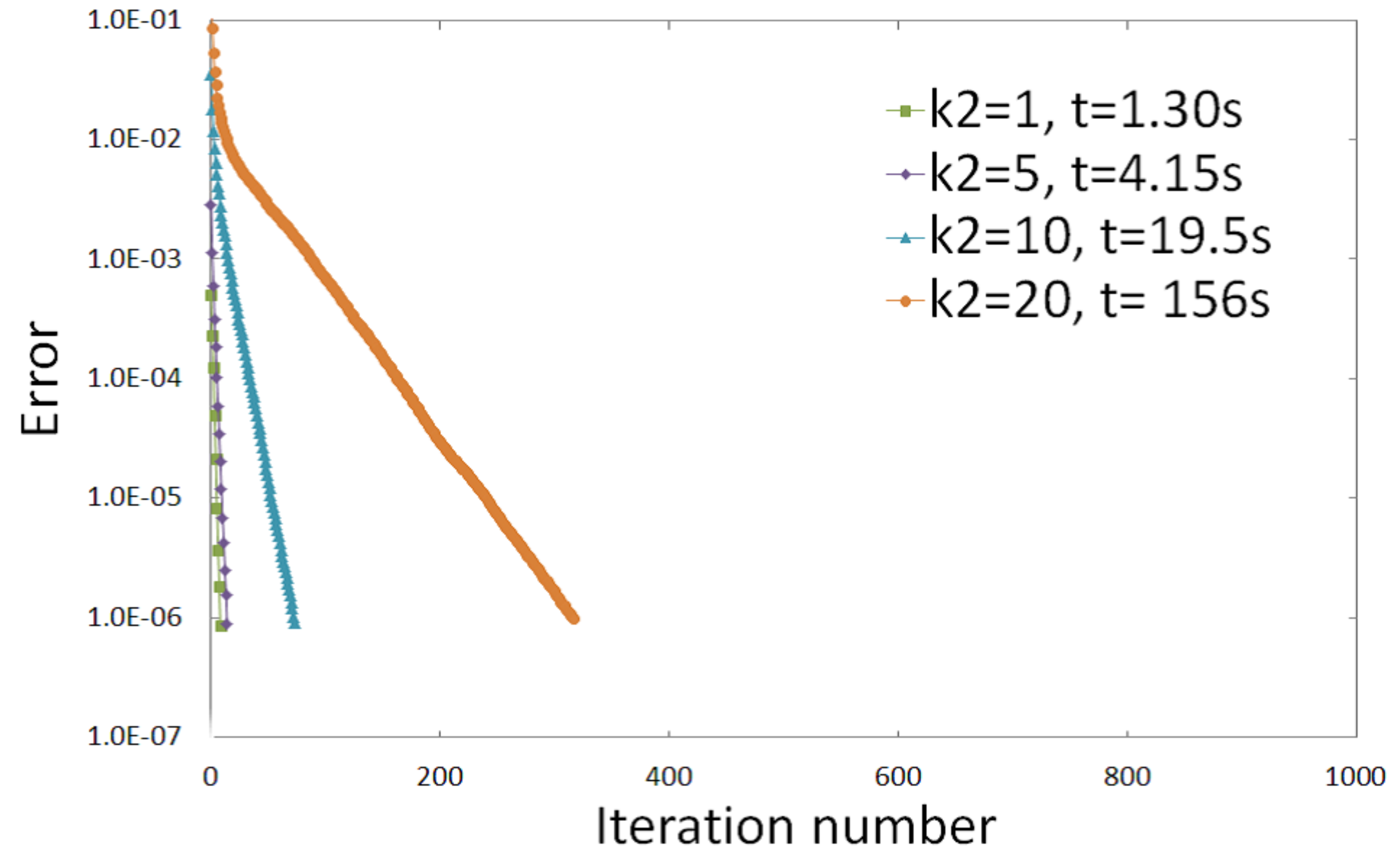}
\subcaption{}
\end{center}
\end{minipage}
\begin{minipage}{0.5\linewidth}
\begin{center}
\includegraphics[width=70mm]{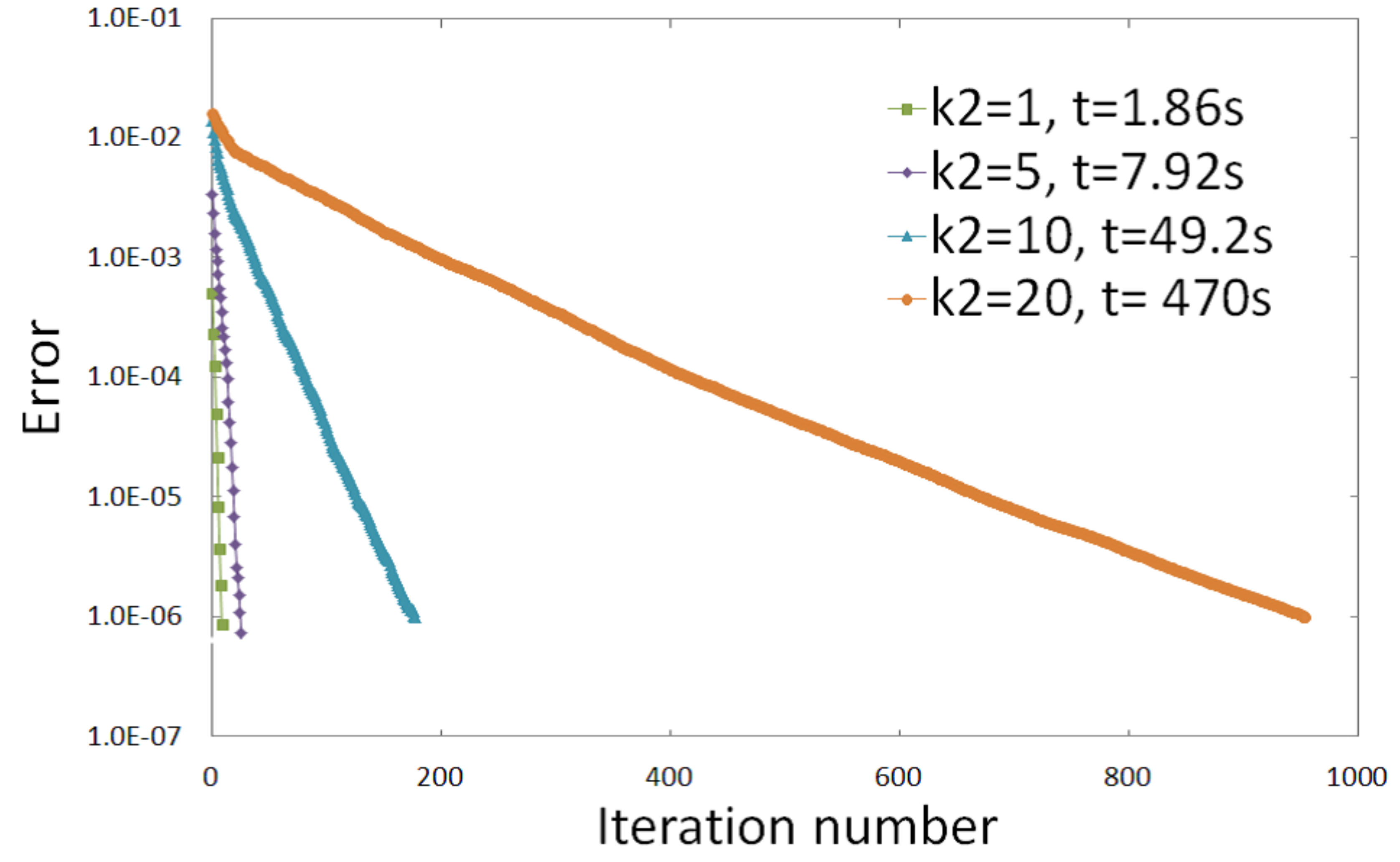}
\subcaption{}
\end{center}
\end{minipage}
\begin{center}
\caption{Convergence behavior of GMRES iteration and the CPU time required
for $200$ inclusions embedded in the central layer, where $k_2$ is allowed to vary
from 1 to 20. In (a), we create a homogeneous background by 
setting $k_1=k_2=k_3$, while in (b), $k_1$ and $k_3$ are fixed at 1, and 
$k_2$ varies.}
\label{fig:NumericalResult2-1}
\end{center}
\end{figure}

In our second example, we consider the same inclusion shape as above, with $k_p=2$,
but with higher contrast materials.
We fix the number of particles to be 200 and the thickness of the middle layer to 
be $d = 12$.  We allow the wavenumber in the middle layer to vary from 
$k_2 = 1$ up to $k_2=20$.  Results are shown in Figs. \ref{fig:NumericalResult2} 
and \ref{fig:NumericalResult2-1} for 6 digits of accuracy.   

In Fig. \ref{fig:NumericalResult2-1}, 
we compare the convergence behavior in an infinite medium (a) vs. a layer medium (b).
Note that the convergence is slower at high constrast and that
this effect is more pronounced in the layerd medium case, where the central 
layer involves strong scattering and reflection.

\subsection{Example 3: scattering from smoothed pentagons}

\begin{figure}[tbp]
\begin{center}
\includegraphics[width=140mm]{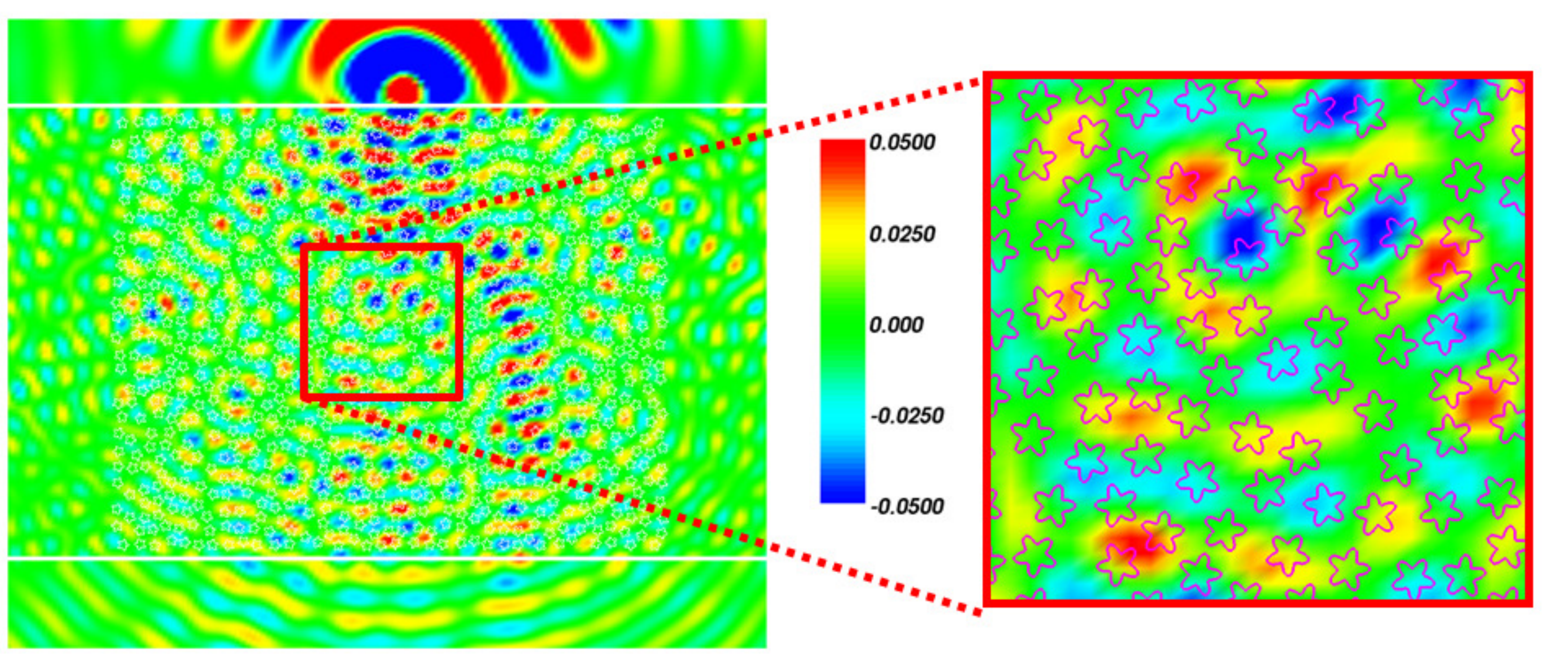}
\caption{Real part of the total field when $1,000$ inclusions are embedded in 
a three-layered medium with $k_1=1.0$, $k_2=3.0$, $k_3=2.0$. 
Each inclusion is a smoothed pentagon, approximately $0.2$ wavelengths in size.}
\label{fig:NumericalResult3}
\end{center}
%
\begin{minipage}{0.5\linewidth}
\begin{center}
\includegraphics[width=70mm]{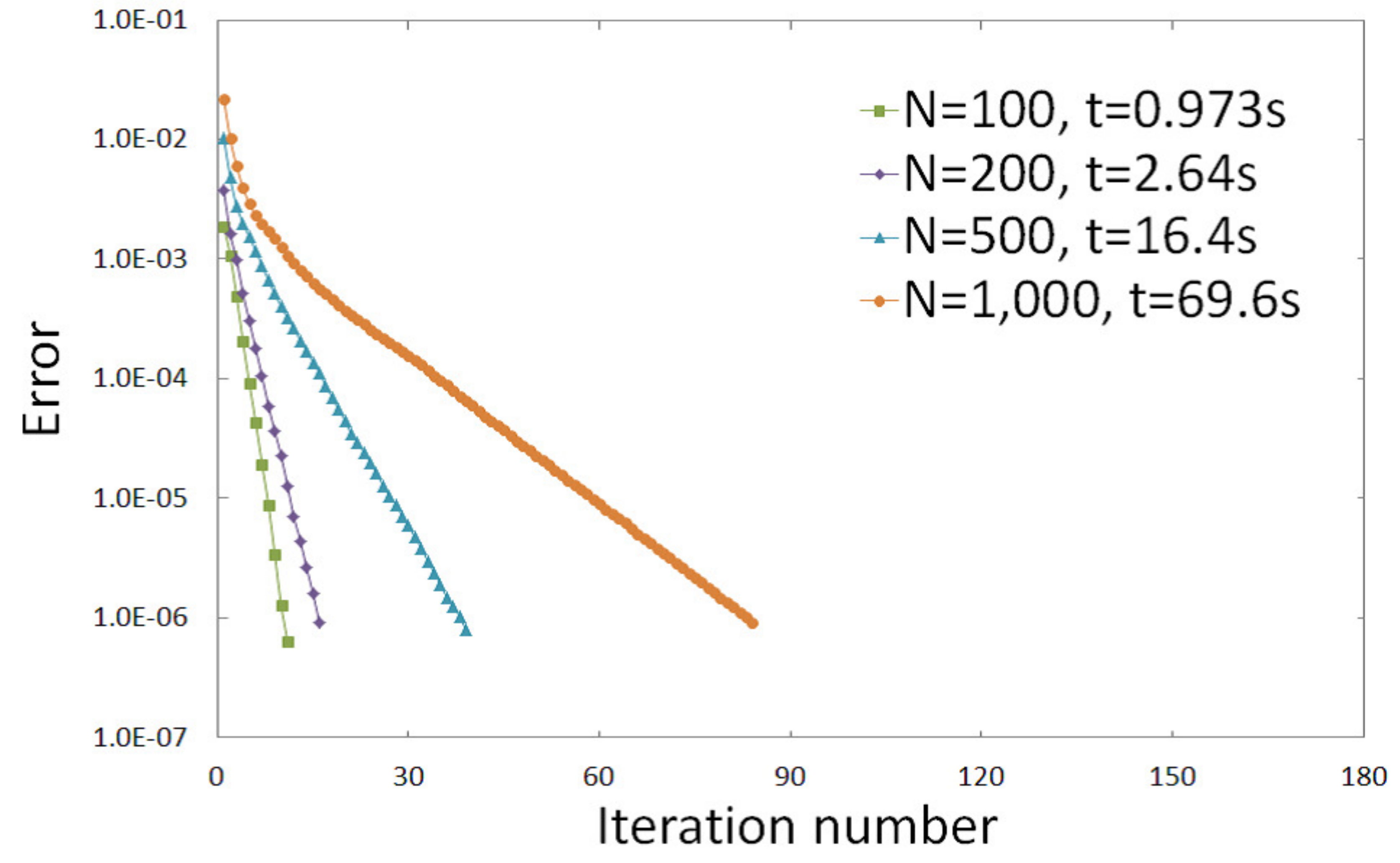}
\subcaption{}
\end{center}
\end{minipage}
\begin{minipage}{0.5\linewidth}
\begin{center}
\includegraphics[width=70mm]{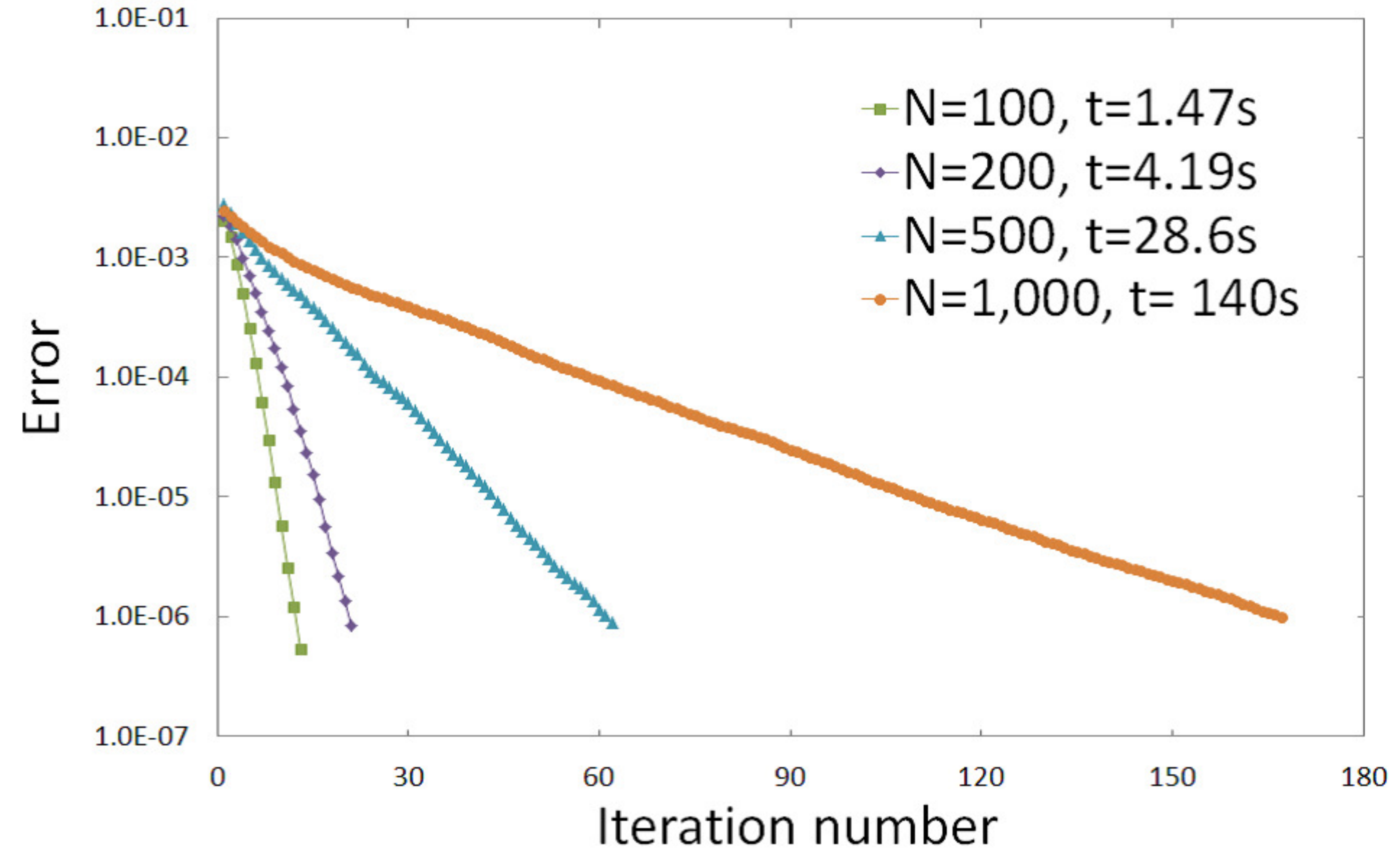}
\subcaption{}
\end{center}
\end{minipage}
\begin{center}
\caption{Convergence behavior of GMRES for a tolerance of $10^{-6}$ and the CPU time required
as the number of inclusions embedded in the central layer varies.
In (a), we create a homogeneous background by 
setting $k_1=k_2=k_3=3.0$, while in (b), $k_1=1.0$, $k_2=3.0$, $k_3=2.0$}
\label{fig:NumericalResult3-1}
\end{center}
\end{figure}

In our last example, we consider the scattering from a different inclusion shape,
setting $a_1=0.3,\ a_2=0.1,\ a_3=5$ in eq. \eqref{par1} with $k_p = 2.0$. 
The inclusions are smoothed pentagons, as shown in Fig. \ref{fig:NumericalResult3}. 
We discretize the boundary of the inclusion using $N=300$ equispaced points and solve eq. \eqref{intg1} 
and \eqref{intg2} to obtain the scattering matrix with $p=10$. 
We consider $M=100$, $200$, $500$ and $1,000$ inclusions. 
For the three-layered medium, we set $\ k_1=1.0,\ k_2=3.0,\ k_3=2.0$.
Results are shown in Figs. \ref{fig:NumericalResult3} and \ref{fig:NumericalResult3-1}.  

Note that in order to obtain 6 digits of accuracy, $69.6$ secs. are required
for $1,000$ inclusions in a homogeneous background, while  
$140$ seconds are required for the three layered medium. 
The time for convergence increases more or less in proportion to $M^2$.

\section{Conclusions}
We have developed a fast algorithm to simulate electromagnetic scattering from 
a microstructured, three-layered material. Our methodology permits inclusions of 
arbitrary shape using a scattering matrix formalism combined with the use of
Sommerfeld integrals to account for the influence of the layered material. 
We have designed efficient procedures to evaluate the Sommerfeld integral at arbitrary 
locations in the layered material using the non-uniform FFT and an effective 
preconditioner that allows the multiple scattering problem to be solved using 
GMRES with a modest number of iterations.
As one would expect from physical considerations, the performance of the method degrades 
when the packing of inclusions is dense and when the contrast is high.
While the method is suitable for parallel implementation, we are also investigating the 
possibility of replacing GMRES iteration with a fast direct solver \cite{GHL2014}.

Extension of the present method to the quasi-periodic case, where the incoming field impinges
on a periodic microstructure will be reported at a later date. 

\section*{Acknowledgements}
This work was supported in part by the Applied
Mathematical Sciences Program of the U.S. Department of Energy
under Contract DEFGO288ER25053 and
by the Office of the Assistant Secretary of Defense for Research and Engineering 
and AFOSR under NSSEFF Program Award FA9550-10-1-0180.

\bibliographystyle{abbrv}
\bibliography{reference}

\end{document}